\documentclass[a4ja,10pt,leqno]{article}

\usepackage{amsmath}
\usepackage{amsfonts}

\addtocounter{equation}{1}
\newcommand{\dbar}{d\!\!\!\lower-0.6ex\hbox{${-}$\!\!}}

\newcommand{\h}{\hbar}

\newcommand{\integ}{\tsize\int\kern-0.7em{\star}\,\,}

\newcommand{\symsum}{\hbox{$\sum
   $\kern-1em\lower-0.3ex\hbox{$\scriptscriptstyle{\bigcirc}$}}\,}

\newcommand{\ott}{\lower-0.4ex\hbox{${\scriptscriptstyle{\otimes}}$}}
\newcommand{\btt}{\lower-0.2ex\hbox{${\scriptscriptstyle{\bullet}}$}}
\newcommand{\ctt}{\lower-0.2ex\hbox{${\scriptscriptstyle{\circ}}$}}
\newcommand{\dtt}{\lower-0.2ex\hbox{${\scriptscriptstyle{\odot}}$}}
\newcommand{\ett}{\lower-0.2ex\hbox{${\scriptscriptstyle{\lozenge}}$}}
\newcommand{\proof}{{\bf{Proof}}\,\,}
\newcommand{\qed}{${}$\hfill$\Box$}
\newtheorem{thm}{Theorem}[section]
\newtheorem{prop}{Proposition}[section]
\newtheorem{lem}{Lemma}[section]

\title{Expressions of algebra elements and transcendental 
noncommutative calculus}

\author{Hideki Omori\\
Department of Mathematics, 
Tokyo University of Science, \\
Noda, Chiba, 278-8510, Japan, \\
\texttt{omori@ma.noda.tus.ac.jp};
\and
Yoshiaki Maeda
\footnote{Partially supported by Grant-in-Aid for 
Scientific Research (\#18204006.), Ministry of Education , Science and 
Culture, Japan.}\\
Department of Mathematics, 
Faculty of Science and Technology, \\
Keio University, Hiyoshi, Yokohama, 
223-8825, Japan, \\
\texttt{maeda@math.keio.ac.jp}; 
\and
Naoya Miyazaki
\footnote{
Partially supported by Grant-in-Aid for 
Scientific Research (\#18540093.), Ministry of Education , Science and 
Culture, Japan.}
\\
Department of Mathematics, Faculty of Economics, \\
Keio University, Hiyoshi, Yokohama, 
223-8521, Japan, \\
\texttt{miyazaki@math.hc.keio.ac.jp}; 
\and
Akira Yoshioka
\footnote{Partially supported by Grant-in-Aid for
  Scientific Research (\#17540096.), Ministry of Education , 
Science and Culture, Japan.}
\\
Department of Mathematics, 
Tokyo University of Science,\\
Kagurazaka, Tokyo, 102-8601, Japan,\\
\texttt{yoshioka@rs.kagu.tus.ac.jp}; 
}

\begin{document}
\maketitle

\begin{center}
\mbox{\parbox{.8\linewidth}{
{\bf Abstract} 
Ideas from deformation quantization are applied to deform the 
expression of elements of an algebra. Extending these ideas to certain 
transcendental elements implies that $\frac{1}{i\h}uv$ in the 
Weyl algebra is naturally viewed as an indeterminate living 
in a discrete set $\mathbb{N}{+}\frac{1}{2}$ {\it or} 
${-}(\mathbb{N}{+}\frac{1}{2})$ . 
This may yield a more mathematical understanding of Dirac's 
positron theory. 
}}  
\end{center}
\par\noindent
{\small{\bf A.M.S Classification (2000)}: 
{Primary 53D55, 53D10; Secondary 46L65}}

\section{Introduction}

Quantum theory is treated algebraically by Weyl algebras, derived from
differential calculus via the correspondence principle. However, since
the algebra is noncommutative, the so-called 
{\it ordering problem} appears.

Orderings are treated in the physics literature of
quantum mechanics (cf. \cite{AW}) as the rules of association from classical 
observables to quantum observables, which are supposed to be 
self-adjoint operators on a Hilbert space. 
Typical orderings are, the normal (standard) ordering, 
the anti-normal (anti-standard) ordering, the Weyl ordering, 
and the Wick ordering in the case of complex variables. 

However, from the mathematical viewpoint, it is better to go back to 
the original understanding of Weyl, which says that orderings are  
procedures of realization of the Weyl algebra $W_\h$. 
Since the Weyl algebra is the universal enveloping algebra of the
Heisenberg Lie algebra, the Poincar{\'e}-Birkhoff-Witt theorem   
shows that this algebra can be viewed as an algebra defined on 
a space of polynomials. As we shall show in \S 1, this indeed 
gives product formulas on the space of polynomials which produce 
algebras isomorphic to $W_\h$. This gives the unique way of 
expressions of elements, and as a result one can treat transcendental 
elements such as exponential functions, which are necessary to solve
differential equations (cf. \S\ref{transcendtal}). 

However, we encounter several anomalous phenomena, such as  
elements with two different inverses 
(cf. \S\ref{anomaloux}) and elements which must be 
treated as double valued 
(cf. \cite{OMMY5},\cite{OMMY6}).  

In this note, we treat the phenomenon which shows that 
$\frac{1}{i\h}uv$ should be viewed as an indeterminate living 
in the set ${\mathbb N}{+}\frac{1}{2}$ {\it or} 
${-}({\mathbb N}{+}\frac{1}{2})$. We reach this interpretation in 
two different ways, by analytic continuation of inverses of 
$z{+}\frac{1}{i\h}uv$, and by defining star 
gamma functions using various ordering
expressions. 

The main point is that we do not use operator theory, but 
instead various ordering expressions, under the leading principle 
that a physical object should be free from 
ordering expressions ({\bf the ordering free principle}), just as
a geometrical object uis free of the local coordinate 
expressions. 

Since similar discrete pictures of elements is familiar in quantum
observables, treated as a self-adjoint operator, our observation 
gives for their justification for the operator theoretic formalism of quantum
theory.   

However, in this note we restrict our ordering expressions 
to a particular subset to avoid the multi-valued expressions. 
In some cases, we should be more careful about the convergence 
of integrals and the continuity of the product, 
so the detailed computations and the proof of continuity 
of the products will appear elsewhere.

\section{$K$-ordering expressions for algebra elements} 

We introcuce a method to realize the Weyl algebra via 
a family of expressions. This leads to 
a transcendental calculus in the Weyl algebra. 

\subsection{Fundamental product formulas and intertwiners}

Let ${\mathfrak S}_{\mathbb C}(n)$ and  
${\mathfrak A}_{\mathbb C}(n)$ be the spaces of complex symmetric 
matrices and skew-symmetric matrices respectively, and 
${\mathfrak M}_{\mathbb C}(n){=}
{\mathfrak S}_{\mathbb C}(n)\oplus{\mathfrak A}_{\mathbb C}(n)$.
For an arbitrary fixed $n{\times}n$-complex matrix 
$\Lambda{\in }{\mathfrak M}_{\mathbb C}(n)$, we 
define a product ${*}_{_{\Lambda}}$ on the space of polynomials   
${\mathbb C}[\pmb u]$ by the formula 
\begin{equation}
 \label{eq:prodgen}
 f*_{_{\Lambda}}g=fe^{\frac{i\h}{2}
(\sum\overleftarrow{\partial_{u_i}}
{\Lambda}{}^{ij}\overrightarrow{\partial_{u_j}})}g
=\sum_{k}\frac{(i\h)^k}{k!2^k}
{\Lambda}^{i_1j_1}\!{\cdots}{\Lambda}^{i_kj_k}
\partial_{u_{i_1}}\!{\cdots}\partial_{u_{i_k}}f\,\,
\partial_{u_{j_1}}\!{\cdots}\partial_{u_{j_k}}g.   
\end{equation}
It is known and not hard to prove that 
$({\mathbb C}[\pmb u],*_{_{\Lambda}})$ is an associative algebra. 

\smallskip
\noindent
(a) The algebraic structure of $({\mathbb C}[\pmb u],*_{_{\Lambda}})$ 
   is determined by the skew-symmetric part of $\Lambda$ (in fact, by
   its conjugacy class $A\to {}^tGAG$).  

\smallskip
\noindent
(b) In particular, if $\Lambda$ is a symmetric matrix, 
  $({\mathbb C}[\pmb u],*_{_{\Lambda}})$ is isomorphic to the usual 
  polynomial algebra. 

\smallskip

Set $\Lambda{=}K{+}J$, $K{\in}{\mathfrak S}_{\mathbb C}(n)$,
$J{\in} {\mathfrak A}_{\mathbb C}(n)$.  
Changing $K$ for a fixed $J$ will be called a {\it deformation} of 
expression of elements, as the algebra remains in the same 
isomorphism class.  

\medskip
\noindent
{\bf Example of computations}:\\
$$
u_i{*_{_{\Lambda}}}u_j
{=}u_iu_j{+}\frac{i\h}{2}{\Lambda}^{ij}, \quad 
u_i{*_{_{\Lambda}}}u_j{*_{_{\Lambda}}}u_k{=}
u_iu_ju_k+\frac{i\h}{2}({\Lambda}^{ij}u_k{+}
{\Lambda}^{ik}u_j{+}{\Lambda}^{jk}u_i).
$$

\medskip
By computing the ${*}_{_{\Lambda}}$-product using the product formula 
\eqref{eq:prodgen}, every element of the algebra has a unique 
expression as a standard polynomial.  We view these expressions 
of an element of algebra as analogous to the ``local coordinate expression'' 
of a function on a manifold. Thus, changing $K$ corresponds 
to a local coordinate transformation on a manifold. 
In this context, we call the product formula \eqref{eq:prodgen} 
the  $K$-{\it ordering expression} by ignoring the 
fixed skew part $J$. For 
$K{=}0, 
\begin{bmatrix}
0&I_m\\
I_m&0\\
\end{bmatrix}, 
\begin{bmatrix}
0&-I_m\\
-I_m&0\\
\end{bmatrix}
$, the $K$-ordering expression is called respectively the Weyl ordering, the normal 
ordering and the anti-normal ordering expressions. 
The intertwiner between a $K$-ordering expression and a 
$K'$-ordering expression,  which we view as a  
local coordinate transformation, is given 
in a concrete form :
\begin{prop}
\label{intwn}
For symmetric matrices $K, K'\in{\mathfrak S}_{\mathbb C}(n)$, the 
intertwiner is given by 
\begin{equation}\label{intertwiner}
I_{_K}^{^{K'}}(f)=\exp\Big(\frac{i\h}{4}\sum_{i,j}(K^{'ij}{-}K^{ij})
\partial_{u_i}\partial_{u_j}\Big)f \,\,
(=I_{0}^{^{K'}}(I_{0}^{^{K}})^{-1}(f)), 
\end{equation}
givieng  an isomorphism 
$I_{_K}^{^{K'}}:({\mathbb C}[{\pmb u}]; *_{_{K{+}J}})\rightarrow 
({\mathbb C}[{\pmb u}]; *_{_{K'{+}J}})$ between algebras.
Namely, for any $f,g \in {\mathbb C}[{\pmb u}]:$ 
\begin{equation}\label{intertwiner2}
I_{_K}^{^{K'}}(f*_{_{K{+}J}}g)=
I_{_K}^{^{K'}}(f)*_{_{K'{+}J}}I_{_K}^{^{K'}}(g).
\end{equation}
\end{prop}

In the case $n{=}2m$ and 
$J{=}
\begin{bmatrix}
  0  &{-}I_m\\
  I_m& 0
\end{bmatrix}, $ 
$({\mathbb C}[\pmb u],*_{_{\Lambda}})$ is called the Weyl
algebra, with isomorphism class denoted by $W_{2m}$. In fact, 
if $J$ is non-singular, then $({\mathbb C}[\pmb u],*_{_{\Lambda}})$
is isomorphic to the Weyl algebra.

\subsection{The star exponential function 
$e_*^{t(z{+}s\frac{1}{i\h}u_k)}$}\label{transcendtal}

Using the ordering expression of elements of algebra, we can treat 
elementary transcendental functions.   
The $*$-exponential function $e_*^{tH}$ is defined as the family 
${:}e_*^{tH}{:}_{_{\Lambda}}$ of solutions of the evolution equations  
\begin{equation}
  \label{eq:evolequ}
\frac{d}{dt}f_{t}{=}H{*_{_{\Lambda}}}f_t, \quad f_0{=}1.
\end{equation}
For instance, for every $z{\in}{\mathbb C}$, we have 
\begin{equation}
  \label{eq:expvv}
{:}e_*^{z{+}s\frac{1}{i\h}u_k}{:}_{_{\Lambda}}
{=}e^z{:}e_*^{s\frac{1}{i\h}u_k}{:}_{_{\Lambda}}
{=}e^{z} e^{s^2\frac{1}{4i\h}K^{kk}}e^{s\frac{1}{i\h}u_k}.  
\end{equation}
When we fix the skew part $J$ of $\Lambda$, we often 
abbreviate the notation to ${:}\,\,{:}_{_K}$, ${*}_{_{K}}$ for 
${:}\,\,{:}_{_{K{+}J}}$, ${*}_{_{K{+}J}}$ respectively.

Since the exponential law 
$$
{:}e_*^{(z{+}w){+}(s{+}t)\frac{1}{i\h}u_k}{:}_{_{K}}{=}
{:}e_*^{z{+}s\frac{1}{i\h}u_k}{:}_{_{K}}{*_{_{K}}}
{:}e_*^{w{+}t\frac{1}{i\h}u_k}{:}_{_{K}}
$$
holds for every $K$, it is better to write  
$$
e_*^{(z{+}w){+}(s{+}t)\frac{1}{i\h}u_k}{=}
e_*^{z{+}s\frac{1}{i\h}u_k}{*}
e_*^{w{+}t\frac{1}{i\h}u_k}
$$
by viewing ${:}e_*^{z{+}s\frac{1}{i\h}u_k}{:}_{_{K}}$ as the
$K$-{\it ordering expression} of the (ordering free) exponential element  
$e_*^{z{+}s\frac{1}{i\h}u_k}$. Under this convention, one may write
for instance
${:}u_i{*}u_j{:}_{_K}{=}u_iu_j{+}\frac{i\h}{2}(K{+}J)^{ij}$.

\medskip
We remark that even for the simplest exponential function 
$e_*^{s\frac{1}{i\h}u_k}$, formula \eqref{eq:expvv} gives the 
following (cf. \cite{om6}). 
\begin{prop}
 \label{theta}
If ${\rm{Im}}\,K^{kk}{<}0$, then the $K$-ordering expression 
of $\sum_{n{\in}\mathbb Z}e_*^{2n\frac{1}{i\h}u_k}$ converges, and 
${:}\sum_{n{\in}\mathbb Z}e_*^{2n\frac{1}{i\h}u_k}{:}_{_K}$ is 
precisely the Jacobi theta function $\theta_3(\frac{1}{i\h}u_k)$. 
\end{prop}
This shows that deformations of expressions  
of a fixed algebraic system are interesting in their own right
(cf. \cite{GG}).
However, it should be remarked that 
$\sum_{n=0}^{\infty} e_*^{2n\frac{1}{i\h}u_k}$, and 
${-}\sum_{n={-\infty}}^{-1}e_*^{2n\frac{1}{i\h}u_k}$ each converge 
to inverses of $1{-}e_*^{\frac{1}{i\h}u_k}$. 
This leads to a breakdown of associativity. Such phenomena occur very often in 
a transcendentally extended algebraic system.   

If ${\rm{Im}}\,K^{kk}{<}0$, then the $K$-ordering expression of 
the integral $\int_{\mathbb R}e_*^{t\frac{1}{i\h}u_k}dt$ converges,
and 
\begin{equation}
\label{eq:identity}
e_*^{z\frac{1}{i\h}u_k}{*}
\int_{\mathbb R}e_*^{t\frac{1}{i\h}u_k}dt{=}
\int_{\mathbb R}e_*^{t\frac{1}{i\h}u_k}dt,\quad 
\forall z{\in}{\mathbb C}.  
\end{equation}
However, we have shown in \cite{OMMY7} that 
$\int_{\mathbb R}e_*^{t\frac{1}{i\h}u_k}dt$ is double valued.

\section{Star exponential functions of quadratic forms}

In this note we mainly deal with the Weyl algebra $W_2$ 
over ${\mathbb C}$. Putting $u_1{=}u, u_2{=}v$, we have 
 the commutation relation $[u,v]{=}-i\h$, 
where $[u,v]{=}u{*}v{-}v{*}u$. The product formula 
\eqref{eq:prodgen} with $\Lambda{=}K{+}J$, 
$J{=}
\begin{bmatrix}
  0&-1\\
  1&0
\end{bmatrix}$ realizes $W_2$.  

In what follows, we use the following notations: 
\begin{equation}
\label{eq:basicform}
u{*}v{=}v{*}u{-}i\h, \quad uv{=}\frac{1}{2}(u{*}v{+}v{*}u),\quad 
v{*}u{=}uv{+}\frac{1}{2}i\h.  
\end{equation}

Let 
$K=\begin{bmatrix}
  0&\kappa\\
  \kappa&0
\end{bmatrix}$. 
The product $*_{_\kappa}$ and the ordering expression $:\,\,:_{_\kappa}$ 
stand for $*_{_K}$ and $:\,\,:_{_K}$, respectively. 
Namely, $*_0$ and $*_1$ correspond to the Moyal product 
and the standard product. 
We also denote the intertwiner 
from the $*_{_\kappa}$-product 
to the $*_{\kappa'}$ product by $I_{\kappa}^{\kappa'}$.

Let $H{\!o}l({\mathbb C}^2)$ be the set of holomorphic functions 
$f(u,v)$ on the complex 2-plane 
${\mathbb C}^2$ endowed with the topology of uniform convergence 
on compact subsets. $H{\!o}l({\mathbb C}^2)$ is viewed as 
a Fr{\'e}chet space.   

The following fundamental lemma flollows easily from  
the product formula \eqref{eq:prodgen}.
 
\begin{lem}
For every polynomial $p(u,v)$, left multiplication $p(u,v){*}$
(resp. right multiplication ${*}p(u,v)$ ) is a continuous 
linear mapping of $H{\!o}l({\mathbb C}^2)$ into itself.  
\end{lem}

\subsection{The star exponential function 
$e_*^{t(z{+}\frac{1}{i\h}uv)}$}\label{subsec3.1}
If $f_t=h(uv)$ in \eqref{eq:evolequ}, 
then 
$I_{\kappa}^{\kappa'}(h(uv))$ is also a function
of $uv$. 
From here on, we mainly concern with functions of $uv$ alone. 
We set $\frac{2}{i\h}uv{=}{\pmb u}A{}^t\!{\pmb u}$, 
where ${\pmb u}{=}(u, v)$ and 
$A{=}
\begin{bmatrix}
  0&1\\
  1&0
\end{bmatrix}$. 
 The intertwiner $I_{\kappa}^{\kappa'}$ is given as follows:
\begin{equation}
  \label{eq:intwner}
I_{\kappa}^{\kappa'}(ge^{t\frac{2}{i\h}uv}){=}
g\frac{1}{1{-}t(\kappa'{-}\kappa)}
e^{\frac{t}{1{-}t(\kappa'{-}\kappa)}\frac{2}{i\h}uv}  
\end{equation}

Solving the evolution equation \eqref{eq:evolequ} 
for the exponential function, we see that 
$e_*^{t\frac{1}{i\h}2uv}$ is given by   
\begin{equation}
  \label{eq:expweyl}
{:}e_*^{t\frac{1}{i\h}2uv}{:}_{0}=
\frac{1}{\cosh t}e^{\frac{1}{i\h}2uv\tanh t}     
\end{equation}
in the Weyl ordering expression (cf. \cite{OMMY5}), and by 
\begin{equation}
  \label{eq:expnorm}
{:}e_*^{t\frac{1}{i\h}2uv}{:}_{_I}
=e^{t}e^{\frac{1}{i\h}(e^{2t}{-}1)uv}  
 \end{equation}
in the normal ordering expression (cf.\cite{OM}). 

Since 
${:}e_*^{t\frac{2}{i\h}uv}{:}_{\kappa}
=I_0^{\kappa}(\frac{1}{\cosh t}e^{\frac{1}{i\h}2uv\tanh t})$, 
we see that
\begin{equation}
\label{kappaexp}
{:}e_*^{t\frac{1}{i\h}2uv}{:}_{\kappa}{=}
\frac{2}{(1{-}\kappa)e^t{+}(1{+}\kappa)e^{-t}}
\exp\big(
\frac{e^t{-}e^{-t}}{(1{-}\kappa)e^t{+}(1{+}\kappa)e^{-t}}
\frac{1}{i\h}2uv\big).  
\end{equation}

Let 
$K=\begin{bmatrix}
  0&\kappa\\
  \kappa&\tau 
\end{bmatrix}$. 
The product $*_{(\kappa,\tau)}$ and 
the ordering expression 
$:\,\,:_{(\kappa,\tau)}$ 
stand for $*_{_K}$ and $:\,\,:_{_K}$, respectively.

It is not hard to obtain the $(\kappa,\tau)$-ordering expression:   
\begin{equation}
  \label{eq:kappatau}
{:}e_*^{t\frac{1}{i\h}2uv}{:}_{(\kappa,\tau)}{=}
\frac{2}{\Delta}
\exp
\big((\frac{e^t{-}e^{-t}}{\Delta})^2\tau\frac{1}{i\h}u^2{+}
\frac{e^t{-}e^{-t}}{\Delta}\frac{1}{i\h}2uv
\big),\,\,\Delta{=}(e^t{+}e^{-t}){-}\kappa(e^t{-}e^{-t}),
\end{equation}
where $\Delta{=}(e^t{+}e^{-t}){-}\kappa(e^t{-}e^{-t})$.
The general ordering expression is a little more complicated involving 
the squre root 
in the amplitude.

\medskip
Note that $(1{-}\kappa)e^t{+}(1{+}\kappa)e^{-t}{=}0$ if and only if
$e^{2t}{=}\frac{\kappa{+}1}{\kappa{-}1}$. Hence, 
${:}e_*^{t\frac{1}{i\h}2uv}{:}_{(\kappa,\tau)}$ has a singular point at 
$2t{=}\log\frac{\kappa{+}1}{\kappa{-}1}{+}2\pi i{\mathbb Z}$. 
However, if $\kappa{=}\pm 1$, then  
${:}e_*^{t\frac{1}{i\h}2uv}{:}_{(\pm 1,\tau)}$ are entire functions with 
respect to  $t$.   
In general we have the following: 

\noindent
\parbox{.6\linewidth}{
\begin{lem}
  \label{exception}
If $\kappa{\in}{\mathbb C}{-}\{\kappa{\geq}1\}
{\cup}\{\kappa{\leq}{-1}\}$, 
then the $(\kappa,\tau)$-ordering expression 
${:}e_*^{t\frac{2}{i\h}uv}{:}_{(\kappa,\tau)}$ is real analytic 
and rapidly decreasing with respect to $t{\in}{\mathbb R}$. 
\end{lem}
}\hfill
\mbox{
\begin{picture}(130,40)(0,20)
\thinlines
\put(0,25){\line(1,0){120}}
\put(60,0){\line(0,1){50}}
\put(33,23){$\bullet$}
\put(80,23){$\bullet$}
\put(63,25){$0$}
\thicklines
\put(82,25){\line(1,0){50}}
\put(0,25){\line(1,0){35}}
\end{picture}}

Formula \eqref{eq:kappatau} gives also the following: 
\begin{prop}
 \label{halfpi}
Suppose $\kappa{\not=}0$, $z{\in}{\mathbb C}$. Then the 
$(\kappa,\tau)$-ordering expression \\
${:}\sin_{*}\pi(z{+}\frac{1}{i\h}uv){:}_{(\kappa,\tau)}$ 
is holomorphic in $(z, uv)$, and vanishes 
on $z{\in}{\mathbb Z}{+}\frac{1}{2}$.
\end{prop}

\noindent
\proof   By \eqref{eq:kappatau}, 
${:}e_*^{\pi i\frac{1}{i\h}2uv}{:}_{(\kappa,\tau)}{+}1{=}0$. 
Although the Weyl ordering expression (the case $\kappa{=}0$) of 
$e_*^{\pm\pi i\frac{1}{i\h}uv}$ diverges by \eqref{eq:expweyl}, 
other ordering expressions exist, e.g. (in normal ordering) 
$$
{:}e_*^{\pi i\frac{1}{i\h}uv}{:}_{1}{=}ie^{-\frac{1}{i\h}2uv},\quad 
{:}e_*^{-\pi i\frac{1}{i\h}uv}{:}_{1}{=}-ie^{-\frac{1}{i\h}2uv}.
$$      
Thus, we have 
$$
0{=}e_*^{-\pi i\frac{1}{i\h}uv}{*}(e_*^{\pi i\frac{1}{i\h}2uv}{+}1){=}
e_*^{\pi i\frac{1}{i\h}uv}{+}e_*^{-\pi i\frac{1}{i\h}uv}{=}
2\cos_*(\pi\frac{1}{i\h}uv).
$$
The desired result follows from the the exponential law. ${}$ \hfill \qed

\begin{lem}
  \label{starsin}
If $\sin_*\pi(z{+}\frac{1}{i\h}uv){*}f(uv)$ is defined on some domain
containing $z{=}\frac{1}{2}$, then 
$\sin_*\pi(\frac{1}{2}{+}\frac{1}{i\h}uv){*}f(uv){=}0$. 
\end{lem}

These observations lead to viewing 
$\frac{1}{2}{+}\frac{1}{i\h}uv$ is an indeterminate in 
the set of integers ${\mathbb Z}$, that is, 
$\frac{1}{i\h}v{*}u$ behaves as if it were an indeterminate 
in ${\mathbb Z}$. However, we have to keep in mind the following 
remark: 

\noindent
{\bf Remark 1}\,\,There are two definitions of the product 
$e_*^{z\frac{1}{i\h}uv}{*}f(u,v)$. The first is to define as the real 
analytic solution of 
$$
\frac{d}{dt}f_t{=}\frac{1}{i\h}uv{*}f_t, \quad f_0{=}f(u,v),
$$
if a real analytic solution exist. The second is to define 
$$
e_*^{z\frac{1}{i\h}uv}{*}f(u,v){=}\lim_{n\to\infty}
e_*^{z\frac{1}{i\h}uv}{*}f_n(u,v),\quad \text{if }\,\,
f(u,v){=}\lim_n f_n(u,v), 
$$ 
where $f_n$ are polynomials. 
These two definitions do not agree in general, since 
the multiplication $e_*^{z\frac{1}{i\h}uv}{*}$ is not a continuous 
linear mapping of $H{\!o}l({\mathbb C}^2)$ into itself 
(cf.\eqref{eq:stridnty}).

\bigskip
\subsection{Several estimates}

${}$

We have already known that 
${:}e_*^{t\frac{1}{i\h}uv}{:}_{\kappa}{\in}H{\!o}l({\mathbb C}^2)$ for
every fixed $t$ whenever defined. 
By \eqref{kappaexp}, we see also that if 
$\kappa{\in}{\mathbb C}{-}\{\kappa{\geq}1\}\cup\{\kappa{\leq}{-1}\}$, 
then ${:}e_*^{t\frac{1}{i\h}uv}{:}_{\kappa}$ is rapidly 
decreasing with respect to $t$.   

In this section, we first show that  
$\int_{-\infty}^{\infty}\,\,e_*^{t\frac{1}{i\h}uv}dt
{\in}H{\!o}l({\mathbb C}^2)$
in the Weyl ordering expression.

\medskip
The Weyl ordering expression of $e_*^{t\frac{1}{i\h}uv}$ is 
${:}e_*^{t\frac{1}{i\h}uv}{:}_0{=}
\frac{1}{\cosh\frac{t}{2}}e^{(\tanh\frac{t}{2})\frac{1}{i\h}2uv}.$
Hence 
$$
{:}\int_{\mathbb R}e_*^{t\frac{1}{i\h}uv}dt{:}_0=
\int_{-\infty}^{\infty}\frac{1}{\cosh\frac{t}{2}}
e^{(\tanh\frac{t}{2})\frac{1}{i\h}2uv}dt.
$$ 
By setting $\cos s{=}\tanh\frac{t}{2}$, $-2\sin sds{=}\sin^2 s dt$,
the integral on the right hand side becomes into  
$$
2\int_{-\pi}^{0}e^{(\cos s)\frac{1}{i\h}2uv}ds{=}
\int_{-\pi}^{\pi}e^{(\cos s)\frac{1}{i\h}2uv}ds.
$$
By the Hansen-Bessel formula, we have  
\begin{equation}
  \label{eq:HansenBessel}
{:}\int_{-\infty}^{\infty}e_*^{t\frac{1}{i\h}uv}dt\,\,{:}_0{=}
\sqrt{\frac{\pi}{2}}J_0(\frac{2}{\h}uv),   
\end{equation}
where $J_0$ is Bessel function of eigen value $0$.

Since $g(s)=e^{(\cos s)\frac{1}{i\h}uv}$ is a continuous curve in 
$H{\!o}l({\mathbb C}^2)$, its integral 
\eqref{eq:HansenBessel} on a compact domain belongs to  
$H{\!o}l({\mathbb C}^2)$. 

Applying the intertwiner $I_0^{\kappa}$ for 
\eqref{eq:HansenBessel}, we see that  
${:}\int_{\mathbb R}
e_*^{t\frac{1}{i\h}uv}dt{:}_{\kappa}=
\int_{-\pi}^{\pi}{:}e^{(\cos s)\frac{1}{i\h}2uv}{:}_{\kappa}ds$.
Since 
$$
{:}e^{(\cos s)\frac{1}{i\h}2uv}{:}_{\kappa}{=}
\frac{2}{(1{-}\kappa)
e^{\frac{1}{2}{\cos s}}{+}(1{+}\kappa)e^{\frac{1}{2}{\cos s}}}
\exp\big(\frac{e^{\cos s}{-}1}
{(1{-}\kappa)e^{\cos s}{+}(1{+}\kappa)}\frac{1}{i\h}2uv\big),
$$
we have the following:
\begin{prop}
\label{estimate}
For every 
$\kappa{\in}{\mathbb C}{-}\{\kappa{\geq}1\}{\cup}\{\kappa{\leq}{-1}\}$, 
the $\kappa$-ordering expression of the integral 
${:}\int_{-\infty}^{\infty}\,\,e_*^{t\frac{1}{i\h}uv}dt{:}_{\kappa}$
is contained in the space $H{\!o}l({\mathbb C}^2)$. 
Furthermore, integration by parts gives 
$\frac{d}{d\theta}\int_{-\infty}^{\infty}\,\,
e_*^{e^{i\theta}t\frac{1}{i\h}uv}
e^{i\theta}dt{=}0$ 
whenever defined.  
\end{prop}

The $*$-delta function is defined by the following integral: 
$$
\int_{-\infty}^{\infty}\,\,e_*^{t\frac{1}{i\h}uv}dt{=}
\int_{\mathbb R}e_*^{-it\frac{1}{\h}uv}dt{=}
\delta_*({\frac{1}{\h}uv}). 
$$

\medskip
Note that $\cos s{=}\tanh\frac{t}{2}$ implies 
$t{=}\log\frac{1{+}\cos s}{1{-}\cos s}$. Hence, we have 
\begin{lem}
\label{replace}
If $f(t)$ is a continuous function such that 
$f(\log\frac{1{+}\cos s}{1{-}\cos s})$ is continuous on $[-\pi,0]$,  
then 
$\int_{-\infty}^{\infty}f(t)e_*^{t\frac{1}{i\h}uv}dt$ is in 
$H{\!o}l({\mathbb C}^2)$ in the $\kappa$-ordering expression
such that for every 
$\kappa{\in}{\mathbb C}{-}\{\kappa{\geq}1\}{\cup}\{\kappa{\leq}{-1}\}$. 
\end{lem}

Applying Lemma\,\ref{replace} to the function 
 $f(t){=}e^{-at}$ $(a{>}0)$ and $e^{-{e^t}}$, we have 
${:}\int_{\mathbb R}
e^{-at}e_*^{t\frac{1}{i\h}uv}dt{:}_{\kappa}$ 
and 
${:}\int_{\mathbb R}e^{-e^{t}}e_*^{t\frac{1}{i\h}uv}dt{:}_{\kappa}$ 
are elements of $H{\!o}l({\mathbb C}^2)$.  
We denote the second integral by  
$$
\int_{\mathbb R}e^{-e^{t}}e_*^{t\frac{1}{i\h}uv}dt{=}
\varGamma_*(\frac{1}{i\h}uv)\quad ({\rm{cf. \, \S\ref{stargamma}}}).
$$

Since $v{*}u{=}uv{+}\frac{1}{2}i\h$, \eqref{eq:kappatau} also gives  
the existence of the limit 
\begin{equation}
  \label{eq:vacuum}
  \begin{aligned}
\lim_{t{\to}\infty}{:}e_*^{t\frac{1}{i\h}2v{*}u}{:}_{(\kappa,\tau)}
&{=}\frac{2}{1{-}\kappa}
e^{\frac{1}{i\h}\frac{1}{1{-}\kappa}
(2uv{+}\frac{\tau}{1{-}\kappa}u^2)},\\ 
\lim_{t{\to}-\infty}{:}e_*^{t\frac{1}{i\h}2u{*}v}{:}_{(\kappa,\tau)}
&{=}\frac{2}{1{+}\kappa}
e^{-\frac{1}{i\h}\frac{1}{1{+}\kappa}
(2uv{-}\frac{\tau}{1{+}\kappa}u^2)},\\
\lim_{t{\to}-\infty}{:}e_*^{t\frac{1}{i\h}2v{*}u}{:}_{(\kappa,\tau)}
&{=}0, 
\quad 
\lim_{t{\to}\infty}{:}e_*^{t\frac{1}{i\h}2u{*}v}{:}_{(\kappa,\tau)}
{=}0.
 \end{aligned}
\end{equation}

We call  
$$
\varpi_{00}{=}\lim_{t{\to}-\infty}e_*^{t\frac{1}{i\h}2u{*}v}, 
\quad 
\overline{\varpi}_{00}{=}\lim_{t{\to}\infty}e_*^{t\frac{1}{i\h}2u{*}v}   
$$
{\bf vacuums}. The exponential law gives 
$$
\varpi_{00}{*}_0\varpi_{00}{=}\varpi_{00}, \quad 
\overline{\varpi}_{00}{*}_0\overline{\varpi}_{00}
{=}\overline{\varpi}_{00}. 
$$
However, we easily see 

\begin{thm}
The product $\varpi_{00}{*}_0\overline{\varpi}_{00}$ diverges 
in any ordering expression.
\end{thm}

The existence of the limit \eqref{eq:vacuum} gives also 
$$
u{*}v{*}\varpi_{00}=0=\varpi_{00}{*}u{*}v.
$$
But the ``bumping identity" $v{*}f(u{*}v){=}f(v{*}u){*}v$ 
give the following: 
\begin{lem}
  \label{vacvac3}
$v{*}\varpi_{00}{=}0{=}\varpi_{00}{*}u$.
\end{lem}

\noindent
\proof 
Using the continuity of $v{*}$, we see that 
$v{*}\lim_{t{\to}-\infty}e_*^{t\frac{1}{i\h}2u{*}v}{=}
\lim_{t{\to}-\infty}v{*}e_*^{t\frac{1}{i\h}2u{*}v}$. 
Hence, the bumping identity (proved by the uniqueness of the 
real analytic solution for linear differential equations) gives 
$\lim_{t{\to}-\infty}e_*^{t\frac{1}{i\h}2v{*}u}{*}v{=}0$ by using 
\eqref{eq:vacuum}. \qed.

However, we note that associativity is not  
easily ensured. The following is the simplest condition
which ensures associativity for certain calculations:  
\begin{prop}
  \label{assoclem}
For every polynomial and for every entire function 
$f{\in} H\!ol({\mathbb C}^2)$, the products 
$p{*}f$ and $f{*}p$ are defined as elements of 
$H\!ol({\mathbb C}^2)$, and associativity 
$(f{*}g){*}h{=}f{*}(g{*}h)$ holds whenever two of 
$f, g, h$ are polynomials.
\end{prop}

In general $(f{*}g){*}h{=}f{*}(g{*}h)$ does not hold 
even if $g$ is a polynomial. 

\noindent
{\bf Example 1}\,\,By Lemma \ref{exception}, 
$\frac{1}{i\h}uv$ has two different inverses 
$$
\big(\frac{1}{i\h}uv\big)_{+}^{-1}{=}\int_{-\infty}^0e_*^{t\frac{1}{i\h}uv}dt, \quad 
\big(\frac{1}{i\h}uv\big)_{-}^{-1}{=}{-}\int_0^{\infty}e_*^{t\frac{1}{i\h}uv}dt 
$$
as elements of $H\!ol({\mathbb C}^2)$. Hence, we see the failure of 
associativity :
$$
\Big((\frac{1}{i\h}uv)_{+}^{-1}{*}
(\frac{1}{i\h}uv)\Big){*}(\frac{1}{i\h}uv)_{-}^{-1}
\not=
(\frac{1}{i\h}uv)_{+}^{-1}{*}
\Big((\frac{1}{i\h}uv){*}(\frac{1}{i\h}uv)_{-}^{-1}\Big),
$$    
and indeed $(\frac{1}{i\h}uv)_{+}^{-1}{*}(\frac{1}{i\h}uv)_{-}^{-1}$
diverges in any ordering expression. 
In what follows, we use the notation
\begin{equation}
 \label{eq:stardelta}
\delta_{*}(\frac{1}{\h}uv){=}
(\frac{1}{i\h}uv)_{+}^{-1}{-}(\frac{1}{i\h}uv)_{-}^{-1}.  
\end{equation}

\medskip

In spite of theis general failure of associativity, we have 
another primitive criterion for associativity. 
We remark that if all terms are 
considered as formal power series in $i\h$ in the product formula 
\eqref{eq:prodgen}, then the product is always defined, 
and it is easy to show associativity, as it holds 
for polynomials (cf. \cite{OM} for details). 
Applying these remarks carefully, we give the following:
\begin{lem}
  \label{vacvac25}
$\varpi_{00}{*}(u^p{*}\varpi_{00}){=}0, and\,\,\,
 (\varpi_{00}{*}v^p){*}\varpi_{00}{=}0$.
\end{lem}

\noindent
\proof\,\,
By taking the formal power series
expansion with respect to $i\h$ for $e_*^{su{*}v}$,  
associativity holds, and the following 
computation is permitted by the bumping identity: 
$$
e_*^{su{*}v}{*}(u^{p}{*}e_*^{tu{*}v}){=}(e_*^{su{*}v}{*}u^{p}){*}e_*^{tu{*}v}
{=}
u^p{*}e_*^{(s{+}t)u{*}v{+}i\h ps}.
$$
The right hand side of the above equality is continuous 
in $s, t$. In particular, 
$$
\lim_{t{\to}a}e_*^{su{*}v}{*}(u^{p}{*}e_*^{tu{*}v}){=}
e_*^{su{*}v}{*}\lim_{t{\to}a}(u^{p}{*}e_*^{tu{*}v}).
$$ 
Using the bumping identity, we have 
$$
\begin{aligned}
e_*^{su{*}v}{*}(u^{p}{*}\lim_{t{\to}{-}\infty}e_*^{tu{*}v})
{=}&e_*^{su{*}v}{*}\lim_{t{\to}{-}\infty}u^{p}{*}e_*^{tu{*}v}
{=}\lim_{t{\to}{-}\infty}u^{p}{*}e_*^{(s{+}t)u{*}v+i\h ps}\\
{=}&u^{p}{*}\lim_{t{\to}{-}\infty}e_*^{(s{+}t)u{*}v+i\h ps}
{=}u^{p}e^{i\h ps}{*}\varpi_{00}.
\end{aligned}
$$
It follows that
$$
\varpi_{00}{*}(u^p{*}\varpi_{00}){=}
\lim_{s{\to}{-}\infty}e_*^{s\frac{1}{i\h}u{*}v}
{*}(\lim_{t{\to}{-}\infty}{u^p}{*}e_*^{t\frac{1}{i\h}u{*}v})
{=}
\lim_{s{\to}{-}\infty}u^{p}e^{ps}{*}\varpi_{00}{=}0.
$$  
Similarly, we also have $(\varpi_{00}{*}v^p){*}\varpi_{00}{=}0$. \qed 

\begin{lem}
  \label{vacpolyvac}
For every polynomial $f(u,v){=}\sum a_{ij}u^i{*}v^{j}$,   
$$
\varpi_{00}{*}(f(u,v){*}\varpi_{00}){=}f(0,0)\varpi_{00}
{=}(\varpi_{00}{*}f(u,v)){*}\varpi_{00}.
$$
Consequently, associativity holds for 
$\varpi_{00}{*}p(u,v){*}\varpi_{00}$ for a polynomial $p(u,v)$ .  
\end{lem}

A similar computation gives the following associativity  
$$
(\varpi_{00}{*}v^q){*}(u^p{*}\varpi_{00}){=}
\delta_{p,q}p!(i\h)^p{=}
\varpi_{00}{*}(v^q{*}u^p{*}\varpi_{00}){=}
(\varpi_{00}{*}v^q{*}u^p){*}\varpi_{00}.
$$
Since 
$$
\varpi_{00}{*}v^q{*}u^{p}{*}\varpi_{00}{=}
\delta_{p,q}p!(i\h)^p\varpi_{00}, 
$$
we have the following:
\begin{prop}
  \label{vacvac4}
$\frac{1}{\sqrt{p!q!(i\h)^{p{+}q}}}u^p{*}\varpi_{00}{*}v^q$ is the 
$(p,q)$-matrix element.
\end{prop}

As mentioned in Remark 1 in \S\,\ref{subsec3.1}, 
we have two definitions of 
$e_*^{z\frac{1}{i\h}uv}{*}f(u,v)$. However both definitions give the 
formula    
\begin{equation}
 \label{eq:stridnty}
e_*^{z\frac{1}{i\h}uv}{*}\varpi_{00}{=}e^{-\frac{1}{2}z}{*}\varpi_{00}.
\end{equation}
 
On the other hand, since  
$\frac{1}{i\h}uv{*}\delta_{*}(\frac{1}{\h}uv){=}0$, 
we must set 
$e_*^{t\frac{1}{i\h}uv}{*}\delta_{*}(\frac{1}{\h}uv)
{=}\delta_{*}(\frac{1}{\h}uv)$ as the real analytic solution of 
$\frac{d}{dt}f_t{=}\frac{1}{i\h}uv{*}f_t$. 

However, computing 
$$
\lim_{N\to\infty}e_*^{t\frac{1}{i\h}uv}{*}\int_{-N}^{N}e_*^{s\frac{1}{i\h}uv}ds
=\lim_{N\to\infty}\int_{-N}^{N}e_*^{(t{+}s)\frac{1}{i\h}uv}ds
$$
gives the following:
\begin{equation}
 \label{eq:stridnty00}
e_*^{(x{+}iy)\frac{1}{\h}uv}{*}\delta_{*}(\frac{1}{\h}uv)
{=}e_*^{{iy}\frac{1}{\h}uv}{*}\delta_{*}(\frac{1}{\h}uv).  
\end{equation}
Hence  \eqref{eq:stridnty} is holomorphic with respect to $z$, while 
\eqref{eq:stridnty00} is only continuous, that is, 
there is no real analyticity with respect to $z=x{+}iy$.

\bigskip

\section{Inverses and their analytic continuation}
\label{anomaloux}

Formula \eqref{eq:expvv} and the exponential law give 
in particular   
$$
{:}e_*^{t(z{+}\frac{1}{i\h}v)}{:}_{(\kappa,\tau)}
{=}e^{\frac{1}{4i\h}t^2\tau}e^{t(z{+}\frac{1}{i\h}v)}.  
$$
It follows that if ${\rm{Im}}\,\tau <0$, then 
$e^{\frac{1}{4i\h}t^2\tau}$ is rapidly decreasing in 
$t$ and the integrals
\begin{equation}
 \label{eq:twoinv00}
{:}\int_{-\infty}^0e_*^{t(z{+}\frac{1}{i\h}v)}dt{:}_{(\kappa,\tau)}, 
\quad 
-{:}\int_{0}^{\infty}e_*^{t(z{+}\frac{1}{i\h}v)}dt{:}_{(\kappa,\tau)}. 
\end{equation}
converge. Both integrals are respectively inverses of $z{+}\frac{1}{i\h}v$,
and are denoted $(z{+}\frac{1}{i\h}v)_{+*}^{-1}$, 
$(z{+}\frac{1}{i\h}v)_{-*}^{-1}$, respectively, with the 
subscript $(\kappa,\tau)$ ommitted.  

\begin{prop}
 \label{property00}
If ${\rm{Im}}\,\tau<0$, then the $(\kappa,\tau)$-ordering 
expression of the difference of the two inverses is given by 
$$
{:}(z{+}\frac{1}{i\h}v)_{+*}^{-1}{-}
(z{+}\frac{1}{i\h}v)_{-*}^{-1}{:}_{(\kappa,\tau)}{=}
\int_{-\infty}^{\infty}
e^{\frac{1}{4i\h}t^2\tau}e^{t(z{+}\frac{1}{i\h}v)}dt.
$$
This difference is holomorphic in $z$. 
\end{prop}

\medskip
Similarly, by formula \eqref{eq:expweyl}, we have the convergence 
of the two integrals
\begin{equation}
\label{+inverse}
{:}\int_{-\infty}^0e^{tz}e_*^{t\frac{1}{i\h}uv}dt{:}_{0}
=\int_{-\infty}^0\frac{e^{\frac{1}{2}tz}}{\cosh\frac{1}{2}t}
e^{\frac{1}{i\h}2uv\tanh\frac{1}{2}t}dt,
\quad {\rm {Re}}\,z>{-\frac{1}{2}},     
\end{equation}
\begin{equation}
\label{-inverse}
{:}{-}\int_{0}^{\infty}e^{tz}e_*^{t\frac{1}{i\h}uv}dt{:}_{0}
={-}\int_{0}^{\infty}
\frac{e^{\frac{1}{2}tz}}{\cosh\frac{1}{2}t}
e^{\frac{1}{i\h}2uv\tanh\frac{1}{2}t}dt,
\quad {\rm {Re}}\,z<{\frac{1}{2}}.     
\end{equation}
Both \eqref{+inverse} and \eqref{-inverse} give inverses 
of $z{+}\frac{1}{i\h}uv$. By a similar computation, there are two
inverses for every $(\kappa,\tau)$ such that 
$\kappa{\in}{\mathbb C}{-}\{\kappa{\geq}1\}{\cup}\{\kappa{\leq}{-1}\}$, 
which will be denoted by 
$(z{+}\frac{1}{i\h}uv)_{+*}^{-1}$, 
$(z{+}\frac{1}{i\h}uv)_{-*}^{-1}$. 

The following may be viewed as a Sato hyperfunction:     
\begin{prop}
  \label{property}
If $-\frac{1}{2}<{\rm {Re}}\,z<\frac{1}{2}$, then the difference 
of the two inverses is given by 
\begin{equation}
  \label{eq:zuv}
(z{+}\frac{1}{i\h}uv)_{+*}^{-1}{-}(z{+}\frac{1}{i\h}uv)_{-*}^{-1}{=}
\int_{-\infty}^{\infty}e_*^{t(z{+}\frac{1}{i\h}uv)}dt.  
\end{equation}
Its $(\kappa,\tau)$-ordering expression is holomorphic on this strip. 
\end{prop}

One can see the right hand side more closely.   
For $-\frac{1}{2}<{\rm{Re}}\,z\leq 0$, the change of variables 
$\tanh\frac{1}{2}t{=}\cos s$ from forms the right hand side 
of \eqref{eq:zuv} into 
$$
2\int_{-\pi}^{0}
(\frac{1{+}\cos s}{1{-}\cos s})^{z}
e^{(\cos s)\frac{1}{i\h}2uv}ds.
$$  
For $0\leq {\rm {Re}}\,z{<}\frac{1}{2}$ and for  
${-}\cos s{=}\tanh\frac{t}{2}$, $2\sin sds{=}\sin^2 s dt$,
the right hand side of \eqref{eq:zuv} transforms into  
$$
2\int_{0}^{\pi}
(\frac{1{+}\cos s}{1{-}\cos s})^{-z}e^{(\cos s)\frac{1}{i\h}uv}ds.
$$
Hence, Lemma \ref{replace} gives that 
$\int_{-\infty}^{\infty}e_*^{t(z{+}\frac{1}{i\h}uv)}dt$ is an 
element of $H{\!o}l({\mathbb C}^2)$.

\bigskip

On the other hand, note that a chang of variables gives 
$$
((-z){+}\frac{1}{i\h}uv)_{-*}^{-1}{=}
{-}\int_0^{\infty}e_*^{{-}t(z{-}\frac{1}{i\h}uv)}dt
{=}{-}\int_{-\infty}^{0}e_*^{(z{-}\frac{1}{i\h}uv)}dt.
$$ 
Thus, we see that
\begin{equation}
  \label{eq:notation}
(z{-}\frac{1}{i\h}uv)_{-*}^{-1}{=}{-}((-z){+}\frac{1}{i\h}uv)_{-*}^{-1}.   
\end{equation}
This is holomorphic on the domain ${\rm {Re}}\,z{>}-\frac{1}{2}$, which
is also the holomorphic domain for $(z{+}\frac{1}{i\h}uv)_{-*}^{-1}$. 

All of these results are easily proved for the 
Weyl ordering expression. However, if 
$\kappa{\in}{\mathbb C}{-}\{\kappa{\geq}1\}\cup\{\kappa{\leq}{-1}\}$,
then 
${:}e_*^{t\frac{1}{i\h}uv}{:}_{\kappa}$ is rapidly decreasing in $t$, 
and the same computation gives the following:
\begin{prop}
  \label{twoinv}
For every $z$ such that ${\rm {Re}}\,z >{-}\frac{1}{2}$, 
the two inverses $(z{+}\frac{1}{i\h}uv)_{+*}^{-1}$ and 
$(z{-}\frac{1}{i\h}uv)_{-*}^{-1}$ are defined in the 
$\kappa$-ordering expression  for 
$\kappa{\in}{\mathbb C}{-}\{\kappa{\geq}1\}\cup\{\kappa{\leq}{-1}\}$. 
\end{prop}

Note that 
$(z{+}\frac{1}{i\h}uv)_{+*}^{-1}{*}({-}z{-}\frac{1}{i\h}uv)_{-*}^{-1}$
diverges for any ordering expression. However, 
the standard resolvent formula gives the following: 
\begin{prop}
  \label{property22}
If $z{+}w{\not=}0$, then  
$$
\frac{1}{z{+}w}
\Big((z{+}\frac{1}{i\h}uv)_{+*}^{-1}{+}(w{-}\frac{1}{i\h}uv)_{-*}^{-1}\Big)
$$
is an inverse of $(z{+}\frac{1}{i\h}uv){*}(w{-}\frac{1}{i\h}uv)$. 
In particular, for every positive integer $n$, and 
for every complex number $z$ such that 
${\rm {Re}}\,z{>}-(n{+}\frac{1}{2})$,
$$
\frac{1}{2n}\big((1{+}\frac{1}{n}(z{+}\frac{1}{i\h}uv))_{+*}^{-1}
{+}(1{-}\frac{1}{n}(z{+}\frac{1}{i\h}uv))_{-*}^{-1}\big)
$$
is an inverse of $1{-}\frac{1}{n^2}(z{+}\frac{1}{i\h}uv)_*^2$ 
in the $\kappa$-ordering expression for  
$\kappa{\in}{\mathbb C}{-}\{\kappa{\geq}1\}{\cup}\{\kappa{\leq}{-1}\}$. 
\end{prop}

\bigskip

\subsection{Analytic continuation of inverses}

${}$

Recall that $(z{\pm}\frac{1}{i\h}uv)_{\pm*}^{-1}$ are holomorphic 
on the domain ${\rm{Re}}\,z >{-}\frac{1}{2}$. It is natural to expect 
that $(z{\pm}\frac{1}{i\h}uv)_{\pm*}^{-1}
{=}C(C(z{\pm}\frac{1}{i\h}uv))_{\pm*}^{-1}$ 
for any non-zero constant $C$. To confirm this, we set
$C{=}e^{i\theta}$ and consider the $\theta$-derivative  
$$ 
e^{i\theta}\int_{-\infty}^{0}e_*^{e^{i\theta}t(z{\pm}\frac{1}{i\h}uv)}dt.
$$
In the $(\kappa,\tau)$-ordering expression, the phase part of the
integrand is bounded  in $t$ and the amplitude is given by 
$$
\frac{2e^{i\theta}tz}{(1{-}\kappa)
e^{e^{i\theta}t/2}+(1{+}\kappa)e^{-e^{i\theta}t/2}},
\quad \kappa{\not=}1.
$$
Hence, the integral converges whenever 
${\rm{Re}}\,e^{i\theta}(z{\pm}\frac{1}{2})>0$, and by integration by
parts this convergence does not depend on $\theta$. It follows that   
$(z{\pm}\frac{1}{i\h}uv)_{\pm*}^{-1}$ are holomorphic 
on the domain ${\mathbb C}{-}\{t; {-}\infty{<}t{<}{-}\frac{1}{2}\}$. 
 
Next, it is natural to expect that the bumping identity 
$(uv){*}v{=}v{*}(uv{-}i\h)$ gives the following ``sliding identities"
$$
v_{+}^{-1}{*}(z{+}\frac{1}{i\h}uv)_{+*}^{-1}{*}v{=}
(z{-}1{+}\frac{1}{i\h}uv)_{+*}^{-1}, \quad 
v_{+}^{-1}{*}(z{-}\frac{1}{i\h}uv)_{-*}^{-1}{*}v{=}
(z{+}1{-}\frac{1}{i\h}uv)_{-*}^{-1}
$$ 
whenever one can use the inverse of $v$ in a suitable ordering
expression. In this section, analytic continuation will be produced via these 
sliding identities.  

\medskip
In this note, we state the sliding identity by using,  
instead of $v^{-1}$, the left inverse $v^{\ctt}$ of 
$v$ given below. 
First of all, we remark that formula \eqref{eq:expweyl} 
also gives  
$$
(u{*}v)_{-*}^{-1}{=}
-\frac{1}{i\h}\int_0^{\infty}e_*^{t\frac{1}{i\h}u{*}v}dt,
\quad
(v{*}u)_{+*}^{-1}{=}
\frac{1}{i\h}\int_{-\infty}^0e_*^{t\frac{1}{i\h}v{*}u}dt.
$$
These gives left/right inverses of $u, v$ 
$$
v^{\ctt}{=}u{*}(v{*}u)_{+*}^{-1},\quad  
u^{\btt}{=}v{*}(u{*}v)_{-*}^{-1},
$$
for it is easy to see that
$$
v{*}v^{\ctt}{=}1,\quad v^{\ctt}{*}v{=}1{-}\varpi_{00}, \quad
u{*}u^{\btt}{=}1,\quad u^{\btt}{*}u{=}1{-}\varpi_{00}.
$$
The bumping identity gives 
$$
v{*}(z{+}\frac{1}{i\h}uv){*}v^{\ctt}{=}z{+}1{+}\frac{1}{i\h}uv, \quad 
v^{\ctt}{*}(z{+}\frac{1}{i\h}uv){*}v{=}
(1{-}\varpi_{00}){*}(z{-}1{+}\frac{1}{i\h}uv).
$$
The successive use of the bumping identity gives  
the following useful formula:
\begin{equation}
  \label{eq:powerinv}
(u{*}(v{*}u)_{+*}^{-1})^n{*}\varpi_{00}
{=}\frac{1}{n!}(\frac{1}{i\h}u)^n{*}\varpi_{00}.
\end{equation}

\medskip
Using $v^{\ctt}$ instead of $v^{-1}$, 
we can give the analytic continuation of inverses. 
However, we have to be careful about the continuity of the 
$*$-product. We compute 
\medskip
$$
\begin{aligned}
v^{\ctt}{*}(z{+}\frac{1}{i\h}uv)_{+*}^{-1}{=}&
u*\int_{{-}\infty}^{0}e_*^{t(\frac{1}{i\h}uv{+}\frac{1}{2})}dt*\!
\int_{{-}\infty}^{0}e_*^{s(z{+}\frac{1}{i\h}uv)}ds\\   
{=}&
u*\int_{{-}\infty}^{0}\int_{{-}\infty}^{0}
e_*^{t(\frac{1}{i\h}uv{+}\frac{1}{2})}{*}
e_*^{s(z{+}\frac{1}{i\h}uv)}dtds\\
{=}&
\int_{{-}\infty}^{0}\int_{{-}\infty}^{0}e^{t\frac{1}{2}{+}sz}
u{*}e_*^{(t{+}s)\frac{1}{i\h}uv}dtds\\
{=}&
\int_{{-}\infty}^{0}\int_{{-}\infty}^{0}e^{t\frac{1}{2}{+}sz-(t{+}s)}
e_*^{(t{+}s)\frac{1}{i\h}uv}{*}u dtds.
\end{aligned}
$$
Hence, we have the identity whenever both sides are defined: 
$$
\begin{aligned}
(v^{\ctt}{*}(z{+}\frac{1}{i\h}uv)_{+*}^{-1}){*}v{=}&
\int_{{-}\infty}^{0}\int_{{-}\infty}^{0}e^{-t\frac{1}{2}{+}s(z{-}1)}
e_*^{(t{+}s)\frac{1}{i\h}uv}{*}(u{*}v)dtds\\
{=}&
\int_{{-}\infty}^{0}(u{*}v){*}e_*^{t\frac{1}{i\h}u{*}v}dt*\!
\int_{{-}\infty}^{0}e_*^{s(z{-}1{+}\frac{1}{i\h}uv)}ds\\
{=}&(1{-}\varpi_{00}){*}(z{-}1{+}\frac{1}{i\h}uv)_{+*}^{-1}.
\end{aligned}
$$
Remarking that 
$$
\varpi_{00}{*}(z{-}1{+}\frac{1}{i\h}uv)_{+*}^{-1}{=}
(z{-}\frac{1}{2})^{-1}\varpi_{00}, 
$$  
whenever $(z{-}1{+}\frac{1}{i\h}uv)_{+*}^{-1}$ is defined, we have   
\begin{equation}
  \label{eq:inverse}
\big(v^{\ctt}{*}(z{+}\frac{1}{i\h}uv)_{+*}^{-1}\big){*}
v{+}(z{-}\frac{1}{2})^{-1}\varpi_{00}
=\big(z{-}1{+}\frac{1}{i\h}uv\big)_{*+}^{-1}.  
\end{equation}
Since $(z{-}\frac{1}{2})^{-1}\varpi_{00}$ is always defined, 
we see that \eqref{eq:inverse} gives the formula for analytic 
continuation. Using this , we have the following 
(see \cite{om5} and \cite{OM} for more details):
\begin{thm}
  \label{contnuation}
The inverses $(z{+}\frac{1}{i\h}uv)_{+*}^{-1}$,
$(z{-}\frac{1}{i\h}uv)_{-*}^{-1}$ extend to holomorphic functions in
$z$ on ${\mathbb C}{-}\{-({\mathbb N}{+}\frac{1}{2})\}$. 
In particular, 
$(z^2{-}(\frac{1}{i\h}uv)^2)_{{\pm}*}^{-1}$ extend to holomorphic functions
of $z$ on this domain. 
\end{thm}

The product 
$(z{+}\frac{1}{i\h}uv)_{+*}^{-1}{*}(w{+}\frac{1}{i\h}uv)_{+*}^{-1}$ is 
naturally defined, but the formula in Theorem \ref{contnuation}
 looks strange at the first 
glance, because 
$z{+}\frac{1}{i\h}uv$ is not zero at $z{=}n{+}\frac{1}{2}$ and  
$(z{+}\frac{1}{i\h}uv)_{+*}^{-1}$ is singular at
$z{=}n{+}\frac{1}{2}$, but 
$(z{+}\frac{1}{i\h}uv){*}(z{+}\frac{1}{i\h}uv)_{+*}^{-1}{=}1$ for 
$z{\not\in}{-}(\mathbb N{+}\frac{1}{2})$.

Note that 
$$
\int_{-\infty}^{0}(z{+}\frac{1}{i\h}uv){*}e_*^{t(z{+}\frac{1}{i\h}uv)}dt
{=}
\left\{
\begin{matrix} 1 & {\rm{Re}}\,z>{-}\frac{1}{2}\\
1{-}\varpi_{00}& z{=}{-}\frac{1}{2}
\end{matrix}\right., 
$$
$$
\int_{-\infty}^{0}
(z{-}\frac{1}{i\h}uv){*}e_*^{t(z{-}\frac{1}{i\h}uv)}dt
{=}
\left\{
\begin{matrix} 1 &{\rm{Re}}\,z>{-}\frac{1}{2}\\
1{-}\overline{\varpi}_{00}& z{=}{-}\frac{1}{2}
\end{matrix}\right.. 
$$
As suggested by these formulas, we extend the definition of 
the $*$-product as follows: 
For every polynomial $p(u,v)$ or $p(u,v){=}e_*^{s\frac{1}{i\h}uv}$, 
\begin{equation}
  \label{eq:gendef}
p(u,v){*}
(z{\pm}\frac{1}{i\h}uv)_{+*}^{-1}{=}
\lim_{N\to\infty}p(u,v)*\!
\int_{-N}^{0}e_*^{t(z{\pm}\frac{1}{i\h}uv)}dt.  
\end{equation}
Hence we have the formula  
\begin{equation}
\label{strange}
(z{+}\frac{1}{i\h}uv){*}(z{+}\frac{1}{i\h}uv)_{+*}^{-1}{=}
\left\{
\begin{matrix} 1 &{\rm{Re}}\,z>{-}\frac{1}{2}\\
1{-}{\varpi}_{00}& z{=}{-}\frac{1}{2}
\end{matrix}\right.. 
\end{equation}
Considering 
$(v^{\ctt})^n{*}(z{+}\frac{1}{i\h}uv)
{*}(z{+}\frac{1}{i\h}uv)_{+*}^{-1}{*}v^n
{=}(v^{\ctt})^n{*}(z{+}\frac{1}{i\h}uv){*}v^n{*}(v^{\ctt})^n
{*}(z{+}\frac{1}{i\h}uv)_{+*}^{-1}{*}v^n$ 
and 
using the formula \eqref{eq:powerinv}, we have the following:
\begin{thm}
 \label{strange}
If we use definition \eqref{eq:gendef} for the $*$-product, then  
\begin{equation}
  \label{eq:remove1}
(z{+}\frac{1}{i\h}uv){*}(z{+}\frac{1}{i\h}uv)_{+*}^{-1}{=}
\left\{
\begin{matrix} 1 & z{\not\in}{-}(\mathbb N{+}\frac{1}{2})\\
1{-}\frac{1}{n!}(\frac{1}{i\h}u)^n{*}\varpi_{00}{*}v^n& 
z{=}{-}(n{+}\frac{1}{2})
\end{matrix}\right.,  
\end{equation}
\begin{equation}
  \label{eq:remove2}
(z{-}\frac{1}{i\h}uv){*}(z{-}\frac{1}{i\h}uv)_{-*}^{-1}{=}
\left\{
\begin{matrix} 1 & z{\not\in}{-}(\mathbb N{+}\frac{1}{2})\\
1{-}\frac{1}{n!}(\frac{1}{i\h}v)^n{*}\overline{\varpi}_{00}{*}u^n&
 z{=}{-}(n{+}\frac{1}{2})
\end{matrix}\right..  
\end{equation}
\end{thm}
Although $z{=}{-}(n{+}\frac{1}{2})$ are all removable singularities 
for \eqref{eq:remove1} and \eqref{eq:remove2}
as a function of $z$, it is better to retain these singular points. 
 
These formulas give in particular for every fixed positive 
integer $m$ 
\begin{equation}
  \label{eq:particular}
(1{+}\frac{1}{m}(z{+}\frac{1}{i\h}uv)){*}
(1{+}\frac{1}{m}(z{+}\frac{1}{i\h}uv))_{+*}^{-1}{=}
\left\{
\begin{matrix} 1 & z{\not\in}{-}(\mathbb N{+}m{+}\frac{1}{2})\\
1{-}\frac{1}{k!}(\frac{1}{i\h}u)^{k}{*}
\varpi_{00}{*}v^{k}& z{=}{-}(k{+}m{+}\frac{1}{2})
\end{matrix}\right.   
\end{equation}
for arbitrary $k\in {\mathbb N}$. 
We state the following identity for later use: 
\begin{equation}
  \label{eq:vanish}
  \varpi_{00}{*}v^n{*}({-}n{-}\frac{1}{2}{+}\frac{1}{i\h}uv){=}
\varpi_{00}{*}(\frac{1}{i\h}u{*}v){*}v^n{=}0.
\end{equation}

\section{An infinite product formula}

Recall the classical formula 
$\sin\pi x{=}\pi x\prod_{k{=}1}^{\infty}(1{-}\frac{x^2}{k^2})$.
Rewrite this as follows:
$$
\prod_{k{=}1}^{\infty}(1{-}\frac{x^2}{k^2}){=}
\frac{1}{2i}\int\chi_{[-\pi,\pi]}(t)e^{itx}dt {=}
\lim_{n\to\infty}
\int\prod^n_{k=1}(1{+}\frac{1}{k^2}\partial_t^2)\delta(t)e^{itx}dt,
$$
where $\chi_{[-\pi,\pi]}(t)$ is the characteristic function of 
the interval $[-\pi,\pi]$. It follows that 

$$
\chi_{[-\pi,\pi]}(t){=}
2i\lim_{n\to\infty}\prod^n_{k=1}
(1{+}\frac{1}{k^2}\partial_t^2)\delta(t)
$$ 
in the space of distributions. 

For $\kappa$ such that  
$|\frac{\kappa{+}1}{\kappa{-}1}|{\not=}1$ , so that 
${:}{e_*^{it\frac{1}{i\h}uv}}{:}_{\kappa}$ is not 
singular on $t\in{\mathbb R}$, we compute as follows: 
$$
\int\chi_{[-\pi,\pi]}(t)
{:}{e_*^{it(z{\pm}\frac{1}{i\h}uv)}}{:}_{\kappa}dt{=}
\int\chi_{[-\pi,\pi]}(t)e^{itz}
{:}{e_*^{{\pm}it\frac{1}{i\h}uv)}}{:}_{\kappa}
dt.
$$
Fixing a cut-off function $\psi(t)$ of compact support such that
$\psi{=}1$ on $[-\pi,\pi]$, we see that 
$$
\int\chi_{[-\pi,\pi]}(t)
{:}{e_*^{t(z{\pm}\frac{1}{i\h}uv)}}{:}_{\kappa}dt{=}
2i\lim_{n{\to}\infty}
\int\prod^n_{k=1}(1{+}\frac{1}{k^2}\partial_t^2)\delta(t)
\psi(t)e^{tz}{:}{e_*^{\pm it\frac{1}{i\h}uv}}{:}_{\kappa}dt.
$$
Integration by parts gives 
$$
\lim_{n{\to}\infty}\int\delta(t)
\prod^n_{k=1}(1{+}\frac{1}{k^2}\partial_t^2)
\psi(t)e^{tz}{:}{e_*^{\pm it\frac{1}{i\h}uv}}{:}_{\kappa}dt{=}
\lim_{n{\to}\infty}\prod^n_{k=1}
{:}(1{+}\frac{1}{k^2}\partial_t^2)
{e_*^{t(z{\pm}\frac{1}{i\h}uv)}}{:}_{\kappa}.
$$
Hence we have in the $\kappa$-ordering expression that 
$$
\int\chi_{[-\pi,\pi]}(t){e_*^{it(z{\pm}\frac{1}{i\h}uv)}}dt{=}
2i\lim_{n{\to}\infty}
\prod^n_{k=1}(1{-}\frac{1}{k^2}(z{\pm}\frac{1}{i\h}uv)^2)_*.
$$
Noting that 
$$
\sin_*\pi(z{\pm}\frac{1}{i\h}uv){=}
\pi(z{\pm}\frac{1}{i\h}uv){*}
\int\chi_{[-\pi,\pi]}(t){e_*^{it(z{\pm}\frac{1}{i\h}uv)}}dt
\in {H\!ol}({\mathbb C}^2),
$$
we have 
\begin{equation}
  \label{eq:sin*}
\sin_*\pi(z{\pm}\frac{1}{i\h}uv){=}  
\pi(z{\pm}\frac{1}{i\h}uv){*}
\lim_{n{\to}\infty}{*}
\prod^n_{k=1}(1{-}\frac{1}{k^2}(z{\pm}\frac{1}{i\h}uv)^2)_*
\end{equation}
in ${H\!ol}({\mathbb C}^2)$.  
In particular, we have 
\begin{prop}\label{infsin*}
In the $\kappa$-ordering expression with 
$|\frac{\kappa{+}1}{\kappa{-}1}|{\not=}1$,
we have 
$$
\sin_*\pi(z{+}\frac{1}{i\h}uv){=}  
\pi(z{+}\frac{1}{i\h}uv){*}
\lim_{n{\to}\infty}\prod^n_{k=1}{*}
(1{-}\frac{1}{k^2}(z{+}\frac{1}{i\h}uv)^2).
$$ 
This is identically zero on the set $z{\in}{\mathbb Z}{+}\frac{1}{2}$.
\end{prop}

The formula in Proposition\,\ref{infsin*} may be rewritten as 
$$
\begin{aligned}
&\sin_*\pi(z{+}\frac{1}{i\h}uv){=}\\
&\pi(z{+}\frac{1}{i\h}uv){*}
\lim_{n{\to}\infty}\prod^n_{k=1}{*}
(1{-}\frac{1}{k}(z{+}\frac{1}{i\h}uv)){*}
e_*^{\frac{1}{k}(z{+}\frac{1}{i\h}uv)}
{*}\prod^n_{k=1}{*}
(1{+}\frac{1}{k}(z{+}\frac{1}{i\h}uv)){*}
e_*^{{-}\frac{1}{k}(z{+}\frac{1}{i\h}uv)}.
\end{aligned}
$$
In \S\ref{stargamma}, we will define a star gamma function via the 
two different inverses mentioned previously and 
give an infinite product formula for the star gamma function. 

\subsection{The product with 
$(z{+}\frac{1}{i\h}uv)_{+*}^{-1}$ and with  
$\big(1{+}\frac{1}{m}(z{+}\frac{1}{i\h}uv)\big)_{+*}^{-1}$}

First we consider the product 
$(z{+}\frac{1}{i\h}uv)_{\pm *}^{-1}{*}\sin_*\pi(z{+}\frac{1}{i\h}uv)$ 
in two different ways. 
One way is by defining: 
\begin{equation}
  \label{eq:prodsin}
  \begin{aligned}
(z{+}\frac{1}{i\h}uv)_{\pm *}^{-1}&{*}\sin_*\pi(z{+}\frac{1}{i\h}uv)\\  
{=}&
\lim_{n\to\infty}(z{+}\frac{1}{i\h}uv)_{\pm *}^{-1}{*}
\Big((z{+}\frac{1}{i\h}uv){*}
\prod^n_{k=1}{*}(1{-}\frac{1}{k^2}(z{+}\frac{1}{i\h}uv)^2)\Big).    
  \end{aligned}
\end{equation}
Since 
$(z{+}\frac{1}{i\h}uv){*}
\prod^n_{k=1}{*}(1{-}\frac{1}{k^2}(z{+}\frac{1}{i\h}uv)^2)$ is a
polynomial, Proposition \ref{assoclem},  
\eqref{eq:remove1} and \eqref{eq:vanish} give     
\begin{equation}
  \label{eq:prodprod}
(z{+}\frac{1}{i\h}uv)_{\pm *}^{-1}{*}\sin_*\pi(z{\pm}\frac{1}{i\h}uv)
{=}\prod_{k=1}^{\infty}{*}(1{-}\frac{1}{k^2}(z{\pm}\frac{1}{i\h}uv)^2).  
\end{equation}

\medskip
The second way is by defining 
\begin{equation}
  \label{eq:firstway11}
(z{+}\frac{1}{i\h}uv)_{\pm *}^{-1}{*}\sin_*\pi(z{+}\frac{1}{i\h}uv)  
{=}
\lim_{N\to\infty}\int_{-N}^{0}
e_*^{t(z{+}\frac{1}{i\h}uv)}{*}\sin_*\pi(z{+}\frac{1}{i\h}uv).
\end{equation}
This may be written as the complex integral 
$$
\frac{1}{2i}\int_{-\infty{+}\pi i}^{0{+}\pi i}
e_*^{t(z{+}\frac{1}{i\h}uv)}dt 
{-}
\frac{1}{2i}\int_{-\infty{-}\pi i}^{0{-}\pi i}
e_*^{t(z{+}\frac{1}{i\h}uv)}dt.
$$
If ${\rm{Re}}z{>}{-}\frac{1}{2}$, then 
adding  
${-}\frac{1}{2}\int_{{-}\pi}^{\pi}e^{it(z{+}\frac{1}{i\h}uv)}dt$ 
to this expressions gives the clockwise contour integral along the
boundary of the domain 
$D{=}\{z{\in}{\mathbb C}; 
{\rm{Re}}\,z{<}0, {-}\pi{<}{\rm{Im}}\,z{<}\pi \}.$ 

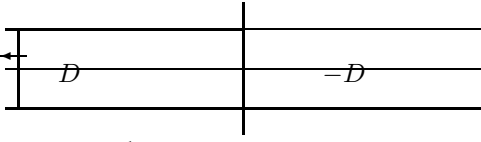
\noindent
\begin{picture}(130,40)(0,15)
\thinlines
\put(0,15){\line(1,0){180}}
\put(90,-10){\line(0,1){50}} 
\put(5,0){\line(0,1){30}}
\put(8,20){\vector(-1,0){10}}
\put(0,0){\line(1,0){180}}
\put(0,30){\line(1,0){180}}
\thicklines
\put(0,0){\line(1,0){90}}
\put(0,30){\line(1,0){90}}
\put(90,0){\line(0,1){30}}
\put(90.02,0){\line(0,1){30}}
\put(20,10){$D$}
\put(120,10){${-}D$}
\end{picture}
\hfill
\parbox{.49\linewidth}{
\begin{lem}
 \label{lem:residue}
${:}e_*^{z\frac{1}{i\h}uv}{:}_{\kappa}$ has at most one singular point
in the domain $D\cup({-}D)$. If ${\rm{Re}}\,\kappa{>}0$, 
then there is no singular point in $D$. 
\end{lem}}

\noindent  
\proof 
${:}e_*^{z\frac{1}{i\h}uv}{:}_{\kappa}{=}
\frac{2}{(1{-}\kappa)e^{\frac{z}{2}}{+}(1{+}\kappa)e^{-\frac{z}{2}}}
\exp{\frac{e^{\frac{z}{2}}{-}e^{-\frac{z}{2}}}
{(1{-}\kappa)e^{\frac{z}{2}}{+}(1{+}\kappa)e^{-\frac{z}{2}}}\frac{2}{i\h}uv}$.
Thus, the singular points are given by     
$(1{-}\kappa)e^{\frac{z}{2}}{+}(1{+}\kappa)e^{-\frac{z}{2}}{=}0$.
This gives 
$e^{z}{=}\frac{\kappa{+}1}{\kappa{-}1}$. If $\kappa{\not=}\pm 1$, then 
$z{=}\log\frac{\kappa{+}1}{\kappa{-}1}{+}2\pi ni$. Thus, 
the domain $D\cup({-}D)$ contains at most one singular point.  

If ${\rm {Re}}\,\kappa{>}0$, then $|\frac{\kappa{+}1}{\kappa{-}1}|>1$
and the singular point (if it exists)  
$z{=}\log\frac{\kappa{+}1}{\kappa{-}1}{+}2\pi ni$ has a 
positive real part. \qed 

\begin{prop}
  \label{prop:residue}
Suppose ${\rm{Re}}\,\kappa{>}0$ and 
$\kappa{\in}{\mathbb C}{-}\{\kappa{\geq}1\}{\cup}\{\kappa{\leq}{-1}\}$.
 Then for $z$ such that ${\rm{Re}}z{>}{-}\frac{1}{2}$, 
we have in the $\kappa$-ordering expression that   
$$
\lim_{N\to\infty}\int_{-N}^{0}
e_*^{t(z{+}\frac{1}{i\h}uv)}{*}\sin_*\pi(z{+}\frac{1}{i\h}uv){=}
\frac{1}{2}\int_{{-}\pi}^{\pi}e_*^{it(z{+}\frac{1}{i\h}uv)}dt.
$$ 
By \eqref{eq:sin*} this integration gives the same result as 
\eqref{eq:prodprod}, namely 
$\prod_{1}^{\infty}{*}(1{-}\frac{1}{k^2}(z{+}\frac{1}{i\h}uv)^2)$.
\end{prop}
 
By the analytic continuation using $v^{\ctt}, v$ as befor, 
we have the following:  
\begin{prop}
  \label{analcontisin}
Suppose ${\rm{Re}}\,\kappa{>}0$ and
$\kappa{\in}{\mathbb C}{-}\{\kappa{\geq}1\}{\cup}\{\kappa{\leq}{-1}\}$. 
Then in the $\kappa$-ordering expression, the product 
$\sin_*\pi(z{+}\frac{1}{i\h}uv){*}(z{+}\frac{1}{i\h}uv)_{+*}^{-1}$ 
is an entire function of $z$. Namely,    
all singularities of $(z{+}\frac{1}{i\h}uv)_{+*}^{-1}$ at
${-}({\mathbb N}{+}\frac{1}{2})$ are cancelled out 
in formulas \eqref{eq:particular} and \eqref{eq:vanish}. 
\end{prop}

By a proof similar to that of Proposition\,\ref{analcontisin}, 
we obtain  
\begin{prop}
\label{inverse}
Suppose ${\rm{Re}}\,\kappa{>}0$, and 
$\kappa{\in}{\mathbb C}{-}\{\kappa{\geq}1\}{\cup}\{\kappa{\leq}{-1}\}$
Then in the $\kappa$-ordering expression,
$$
\sin_{*}\pi(z{-}\frac{1}{i\h}uv){*}(z{-}\frac{1}{i\h}uv)_{-*}^{-1}
$$
is a well defined entire function of $z$.  

In particular, 
$\sin_*\pi(z{+}\frac{1}{i\h}uv){*}(z^2{-}(\frac{1}{i\h}uv)^2)_{{\pm}*}^{-1}$
is a holomorphic function of $z$ in ${\mathbb C}$. 
\end{prop}

\bigskip
Consider next the product 
$(1{+}\frac{1}{m}(z{+}\frac{1}{i\h}uv))_{+*}^{-1}
{*}\sin_*\pi(z{+}\frac{1}{i\h}uv)$.
Since  
$$
(1{+}\frac{1}{m}(z{+}\frac{1}{i\h}uv))_{+*}^{-1}{=}
m(m{+}z{+}\frac{1}{i\h}uv)_{+*}^{-1},
$$
and 
$\sin_*\pi(z{+}m{+}\frac{1}{i\h}uv)
{=}(-1)^m\sin_*\pi(z{+}\frac{1}{i\h}uv)$
by the exponential law, the product formula is 
essentially the same as above.
Hence we see the following:

\begin{prop}\label{prop5.5}
Suppose ${\rm{Re}}\,\kappa{>}0$, and 
$\kappa{\in}{\mathbb C}{-}\{\kappa{\geq}1\}{\cup}\{\kappa{\leq}{-1}\}$. 
Then in the $\kappa$-ordering expression,  the product 
$\sin_*\pi(z{+}\frac{1}{i\h}uv)
{*}(1{+}\frac{1}{m}(z{+}\frac{1}{i\h}uv))_{+*}^{-1}$ 
is an entire function of $z$ with no removable singularity.  
\end{prop}

\noindent
{\bf Remark 2}\,\,
Suppose ${\rm{Re}}\,\kappa{<}0$, and 
$\kappa{\in}{\mathbb C}{-}\{\kappa{\geq}1\}{\cup}\{\kappa{\leq}{-1}\}$. 
Then the residue of 
$e_*^{t(z{+}\frac{1}{i\h}uv)}$ at the singular point 
$t{=}\log\frac{\kappa{+}1}{\kappa{-}1}{+}2\pi ni$
in $D$ gives the difference between the twosides of 
the equality in Proposition\,\ref{analcontisin}. 

This observation shows that 
continuity does not hold for $\kappa$-ordering expressions. 

Lemma \ref{lem:residue} and formula \eqref{kappaexp} show that 
the integral 
$\frac{1}{2\pi i}\int_{\partial D}e_*^{t(z{+}\frac{1}{i\h}uv)}dt$ 
gives the residue at the singular point in $D$. 
This residue will be computed in the last section.

\section{Star gamma functions}\label{stargamma}

We first recall the ordinary gamma function and beta function: 
$$
\varGamma(z){=}\int_0^{\infty}e^{-t}t^{z-1}dt, \quad
B(x,y){=}\int_0^{1}t^{x{-}1}(1{-}t)^{y{-}1}dt. 
$$ 
Substituting $t=e^s$ gives 
$$
\varGamma(z)=\int_{-\infty}^{\infty}e^{-e^{s}}e^{sz}ds,  \quad
B(x,y){=}\int_{-\infty}^{0}e^{sx}(1{-}e^s)^{y{-}1}ds.
$$
The star gamma function and the star beta function may be 
defined by replacing $x$ with $z\pm\frac{uv}{i\h}$:
\begin{equation}
\label{eq:Gamm}
\begin{aligned}[c]
\varGamma_*(z\pm\frac{uv}{i\h})=&
\int_{-\infty}^{\infty}
e^{-e^{\tau}}e_*^{\tau(z\pm\frac{uv}{i\h})}d\tau, \\    
B_*(z\pm\frac{uv}{\h i}, y)= &
\int_{-\infty}^0
e_*^{\tau(z\pm\frac{uv}{i\h})}(1{-}e^\tau)^{y{-}1}d\tau. 
\end{aligned}
\end{equation}
The Weyl ordering expressions of these orderings are 
$$
{:}\varGamma_*(z\pm\frac{uv}{i\h}){:}_{0}
=\int_{-\infty}^{\infty}
\frac{e^{-e^{\tau}{+}z\tau}}{\cosh{\frac{1}{2}\tau}}
e^{\pm\frac{1}{i\h}uv\tanh{\frac{1}{2}\tau}}d\tau,
$$
$$
:B_*(z\pm\frac{uv}{i\h}, y):_0=
\int_{-\infty}^0\frac{(1{-}e^\tau)^{y{-}1}e^{\tau z}}
{\cosh{\frac{1}{2}\tau}}
e^{\pm\frac{1}{i\h}uv\tanh{\frac{1}{2}\tau}}d\tau.
$$
The $\kappa$-ordering expressions  are obtained by 
applying the intertwiner $I_0^{\kappa}$ for  
$\kappa{\in}{\mathbb C}{-}\{\kappa{\geq}1\}{\cup}\{\kappa{\leq}{-1}\}$.
 
\begin{equation}
  \label{eq:gamma}
:\varGamma_*(z\pm\frac{uv}{\h i}):_{\kappa}=
\lim_{N,N'\to\infty}\int_{-N}^{N'}
\frac{e^{-e^{\tau}{+}z\tau}}{\cosh{\frac{1}{2}\tau}}
I_{0}^{\kappa}(e^{\pm\frac{1}{\h i}uv\tanh{\frac{1}{2}\tau}})d\tau.  
\end{equation}
The right hand side converges on a dense open domain of $\kappa$.
\begin{prop}
For every $uv\in {\mathbb C}$, and for every 
$z\in {\mathbb C}$  such
that ${\rm Re}\,z>-\frac{1}{2}$, the right hand side of \eqref{eq:gamma}
converges and is holomorphic with respect to $z$ . However,
$\varGamma_*(-\frac{1}{2}\pm\frac{uv}{\h i})$ is singular. 
\end{prop}

Throughout this section, ordering expressions are always restricted 
to $\kappa{\in}{\mathbb C}{-}\{\kappa{\geq}1\}{\cup}\{\kappa{\leq}{-1}\}$.

\subsection{Analytic continuation of 
$\varGamma_*(z\pm\frac{uv}{\h i})$}

As with the usual gamma function, 
integration by parts gives the identity 
\begin{equation}
 \label{eq:func}
\varGamma_*(z{+}1\pm\frac{uv}{\h i}){=}
(z\pm\frac{uv}{\h i}){*}
\varGamma_*(z\pm\frac{uv}{\h i}).
\end{equation}
Using  
$$
\varGamma_*(z\pm\frac{uv}{\h i}){=}
(z\pm\frac{uv}{\h i})_{\pm*}^{-1}{*}
\varGamma_*(z{+}1\pm\frac{uv}{\h i}),
$$ 
and careful treating continuity insures the following 
\begin{prop}
 \label{conanal}
$\varGamma_*(z\pm\frac{uv}{\h i})$ extends to a holomorphic 
function on 
$z\in {\mathbb C}{-}\{-({\mathbb N}{+}\frac{1}{2})\}$.  
\end{prop}

Since 
$e_*^{\tau(z\pm\frac{uv}{\h i})}{*}\varpi_{00}{=}
(z\pm{\textstyle{\frac{1}{2}}})^{-1}\varpi_{00}$, 
we see the following remarkable feature of 
these star functions 
\begin{equation}
  \label{eq:featue}
  \begin{aligned}[c]
\varGamma_*(z\pm\frac{uv}{\h i})
{*}\varpi_{00}\equiv &
\lim_{N\to\infty}\int_{-N}^Ne^{-e^{\tau}}
e_*^{\tau(z\pm\frac{uv}{\h i})}d\tau 
{*}\varpi_{00}
 =\varGamma(z\pm{\textstyle{\frac{1}{2}}})\varpi_{00} \\
B_*(z\pm\frac{uv}{\h i}, y){*}\varpi_{00}\equiv &
\lim_{N\to\infty}\int_{-N}^N
e_*^{\tau(z\pm\frac{uv}{\h i})}(1{-}e^\tau)^{y{-}1}d\tau{*}\varpi_{00}
=B(z\pm{\textstyle{\frac{1}{2}}}, y)\varpi_{00}    
 \end{aligned}
\end{equation}

\subsection{An infinite product formula}

We see in the same notation as above  
\begin{equation}
  \label{eq:trivial-}
B_*(z{\pm}\frac{uv}{\h i},1)=
\int_{-\infty}^0\!e_*^{\tau(z{\pm}\frac{uv}{\h i})}d\tau
=\big(z{+}\frac{uv}{i\h}\big)_{*\pm}^{-1}, 
\quad {\rm{Re}}\,z>-\frac{1}{2}.  
\end{equation}
We now compute 
$$
{\varGamma}_*(z\pm\frac{uv}{\h i}){\varGamma}(y)=
\iint_{{\mathbb R}^2}e_*^{\tau(z\pm\frac{uv}{\h i})}e^{\sigma y}
e^{-(e^{\tau}{+}e^{\sigma})}d\tau d\sigma.
$$
We change variables by setting 
$$
\tau=t{+}s, \quad e^\sigma=e^{t}(1{-}e^{s}),
\quad{\text{where}}\,\,-\infty<t<\infty, \,\,-\infty<s<0.
$$
Since $e^{\tau}{+}e^{\sigma}=e^t$, this gives a diffeomorphism of
${\mathbb R}\times {\mathbb R}_{-}$ onto ${\mathbb R}^2$. 
The Jacobian is given by 
$d\tau d\sigma= \frac{1}{1{-}e^s}dtds$. 
Hence we have the fundamental relation between the gamma 
function and the beta function 
\begin{equation}
  \label{eq:betgmm}
  \begin{aligned}[c]
{\varGamma}_*(z\pm\frac{uv}{\h i}){\varGamma}(y)
=&\int_{-\infty}^{\infty}\int_{-\infty}^0
e_*^{t(y{+}z\pm\frac{uv}{\h i})}e^{-e^t}{*}
e_*^{s(z\pm\frac{uv}{\h i})}
(1{-}e^s)^{y{-}1}dtds\\
=&{\varGamma}_*(y{+}z\pm\frac{uv}{\h i}){*}B_*(z\pm\frac{uv}{\h i},y).
  \end{aligned}
\end{equation}

Integration by parts gives 
$$
(z\pm\frac{uv}{\h i}){*}B_*(z\pm\frac{uv}{\h i}, y{+}1)=
yB_*(1{+}z\pm\frac{uv}{\h i}, y{+}1).
$$
To prove this, note that 
$$
\frac{d}{d\tau}e_*^{\tau(z\pm\frac{uv}{\h i})}=
(z\pm\frac{uv}{\h i}){*}e_*^{\tau(z\pm\frac{uv}{\h i})}, \quad 
\frac{d}{d\tau}e^{-e^\tau}=-e^{\tau}e^{-e^\tau}, 
$$
$$
\lim_{\tau\to\pm\infty}e^{-e^{\tau}{+}z\tau}e_*^{\pm \tau\frac{uv}{\h i}}=0 
\quad {\text{for}}\quad {\rm{Re}}\,z>-\frac{1}{2}. 
$$
Since 
$B_*(z\pm\frac{uv}{\h i}, y{+}1)
=B_*(z\pm\frac{uv}{\h i}, y){-}B(1{+}z\pm\frac{uv}{\h i},y)$, we have 
the functional equation 
\begin{equation}
  \label{eq:funcrel}
B_*(z\pm\frac{uv}{\h i}, y)
=\frac{y{+}z\pm\frac{uv}{\h i}}{y}{*}B_*(z\pm\frac{uv}{\h i}, y{+}1).
\end{equation}
Iterate \eqref{eq:funcrel} to obtain 
$$
B_*(z\pm\frac{uv}{\h i}, y)=
\frac{(y{+}z\pm\frac{uv}{\h i}){*}
(y{+}1{+}z\pm\frac{uv}{\h i})}{y(y{+}1)}
{*}B_*(z\pm\frac{uv}{\h i}, y{+}2).
$$
Using the notation 
$$
(a)_n=a(a{+}1)\cdots(a{+}n{-}1), \quad \{A\}_{*n}=A{*}(A{+}1){*}\cdots{*}(A{+}n{-}1),
$$
we have 
\begin{equation}
  \label{eq:iteite1}
B_*(z\pm\frac{uv}{\h i}, y)=
\frac{\{y{+}z\pm\frac{uv}{\h i}\}_{*n}}{(y)_n}{*}
B_*(z\pm\frac{uv}{\h i}, y{+}n).
\end{equation}
Similarly, integration by parts gives the formula
\begin{equation}
  \label{eq:fcteq}
\varGamma_*(1{+}z\pm\frac{uv}{\h i})=
(z\pm\frac{uv}{\h i}){*}\varGamma_*(z\pm\frac{uv}{\h i}), 
\quad {\text{for}}\quad {\rm{Re}}\,z>-\frac{1}{2}.  
\end{equation}
Iterate \eqref{eq:fcteq} to obtain  
\begin{equation}
  \label{eq:iteite2}
\varGamma_*(n{+}1{+}z\pm\frac{uv}{\h i})=
\varGamma_*(z\pm\frac{uv}{\h i}){*}
\big\{z\pm\frac{uv}{\h i}\big\}_{*n}.
\end{equation}

\begin{lem}
  \label{prod10}
$B_*(z\pm\frac{uv}{\h i}, n{+}1)=
n!\prod_{k=0}^n{*}(k{+}z\pm\frac{uv}{\h i})^{-1}_{\pm *}$.
\end{lem}

\noindent
\proof  The right hand side of the above equality will be denoted by 
$\frac{n!}{\{z\pm\frac{uv}{\h i}\}^{(\pm)}_{*{n{+}1}}}$.

The case $n=0$ is given by  \eqref{eq:trivial-}. Suppose the formula holds
for $n$. For the case $n{+}1$, we see that 
$$
B_*(z\pm\frac{uv}{\h i}, n{+}2)=
\int_{-\infty}^{0}
e_*^{\tau(z\pm\frac{uv}{\h i})}(1{-}e^{\tau})(1{-}e^{\tau})^nd\tau
=\frac{n!}{\{z\pm\frac{uv}{\h i}\}^{(\pm)}_{*{n{+}1}}}
-\frac{n!}{\{1{+}z\pm\frac{uv}{\h i}\}^{(\pm)}_{*{n{+}1}}}. 
$$
It follows that
$$
B_*(z\pm\frac{uv}{\h i}, n{+}2)=
\frac{(n{+}1)!}{\{z\pm\frac{uv}{\h i}\}^{(\pm)}_{*{n{+}2}}}.
$$
${}$ \hfill $\Box$

In this subsection, we give an infinite product formula for the
$*$-gamma function. 
By Lemma\,\ref{prod10}, we see that
$$
\int_{-\infty}^0
e_*^{\tau(z\pm\frac{uv}{\h i})}(1-e^{\tau})^nd\tau =
\frac{n!}{\{z\pm\frac{uv}{\h i}\}^{(\pm)}_{*{n{+}1}}},
\quad {\rm{Re}}\,z>-\frac{1}{2}.
$$
Replacing $e^{\tau}$ by $\frac{1}{n}e^{\tau'}$, namely setting 
$\tau=\tau'-\log n$ in the left hand side,  and 
multiplying both side by $e_*^{(\log n)(z\pm\frac{uv}{\h i})}$,
we have   
\begin{equation}
  \label{eq:form01}
\int_{-\infty}^{\log n}e_*^{\tau'(z\pm\frac{uv}{\h i})}
(1{-}\frac{1}{n}e^{\tau'})^nd\tau' 
=\frac{n!}{\{z\pm\frac{uv}{\h i}\}^{(\pm)}_{*{n{+}1}}}
{*}e_*^{(\log n)(z\pm\frac{uv}{\h i})}.  
\end{equation}
\begin{lem}
\label{keytop}  
The Weyl ordering expression of the left hand side of 
\eqref{eq:form01} converges 
as $n{\to}\infty$ to 
$\int_{-\infty}^{\infty}
e_*^{\tau'(z\pm\frac{uv}{\h  i})}e^{-e^{\tau'}}d\tau'$ 
in $H{\!o}l({\mathbb C}^2)$. 
\end{lem}

\noindent
{\proof}\,Obviously, 
$\lim_{n{\to}\infty}(1{-}\frac{1}{n}e^{\tau'})^n{=}e^{-e^{\tau'}}$ 
uniformly on each compact subset as a function of $\tau'$. 
In the Weyl ordering expression, it is easy to show that 
$$
\lim_{n\to\infty}
\int_{-\infty}^{\log n}
e_*^{\tau'(z\pm\frac{uv}{\h i})}e^{-e^{\tau'}}d\tau'
{=}
\int_{-\infty}^{\infty}e_*^{\tau'(z\pm\frac{uv}{\h i})}e^{-e^{\tau'}}d\tau'
$$
in $H{\!o}l({\mathbb C}^2)$.
Thus it is enough to show that 
$$
\lim_{n\to\infty}\int_{-\infty}^{\log n}
e_*^{\tau'(z\pm\frac{uv}{\h i})}
(e^{-e^{\tau'}}{-}(1{-}\frac{1}{n}e^{\tau'})^n)d\tau'{=}0
$$
in $H{\!o}l({\mathbb C}^2)$. 
This is easy in the Weyl ordering. 
Applying the intertwiner gives the desired result.  \hfill $\Box$

\medskip
The right hand side of \eqref{eq:form01} equals 
$$ 
e_*^{(\log n{-}(1{+}\frac{1}{2}{+}\cdots{+}\frac{1}{n}))
(z\pm\frac{uv}{\h i})} 
{*}(z\pm\frac{uv}{\h i})_{*\pm}^{-1}{*}
\prod_{k=1}^n
\Big(\big(1{+}\frac{z\pm\frac{uv}{\h i}}{k}\big)_{\pm *}^{-1}
  {*}e_*^{{\frac{z\pm\frac{uv}{\h i}}{k}}}\Big),\quad 
{\rm{Re}}\,z>-\frac{1}{2}. 
$$
The left hand side converges,  and 
$\lim_{n{\to}\infty}
e_*^{(\log n{-}(1{+}\frac{1}{2}{+}\cdots{+}\frac{1}{n}))
(z\pm\frac{uv}{\h i})}
=e_*^{-\gamma (z\pm\frac{uv}{\h i})}$ obviously, where  
$\gamma$ is Euler's constant. By the continuity of the
$*$-multiplication $e_*^{s\frac{uv}{\h i}}{*}$, 
we have the convergence in $H{\!o}l({\mathbb C}^2)$ of 
$$
\lim_{n{\to}\infty}
\prod_{k=1}^n{*}
\Big(\big(1{+}\frac{1}{k}(z\pm\frac{uv}{\h i})\big)_{\pm *}^{-1}
 {*}e_*^{\frac{1}{k}(z\pm\frac{uv}{\h i})}\Big).
$$
Hence we have the convergence in $H{\!o}l({\mathbb C}^2)$ of
the infinite product formula  
\begin{equation}
  \label{eq:prodinf}
{\varGamma}_*(z{+}\frac{uv}{\h i})=
e_*^{-\gamma(z{+}\frac{uv}{\h i})}
{*}(z{+}\frac{uv}{\h i})_{*+}^{-1}{*}
\prod_{k=1}^{\infty}{*}
\Big(\big(1{+}\frac{1}{k}(z{+}\frac{1}{i\h}uv)\big)_{+*}^{-1}
  {*}e_*^{\frac{1}{k}(z{+}\frac{1}{i\h}uv)}\Big)
\end{equation}

Fix $m{\in}{\mathbb N}$. Multiplying 
$(1{+}\frac{1}{m}(z{+}\frac{uv}{\h i})
e_*^{-{\frac{1}{m}(z{+}\frac{uv}{\h i})}}$ 
to both side of \eqref{eq:prodinf}  
and using the abbreviated notation 
$$
\prod_{k{\not=}m}(z{=}a)\,\,{=}\,\,
(z{+}\frac{uv}{\h i})_{+*}^{-1}
\prod_{k{\not=}m}{*}
\Big(\big(1{+}\frac{1}{k}(a{+}\frac{1}{i\h}uv)\big)_{+*}^{-1}
  {*}e_*^{\frac{1}{k}(a{+}\frac{1}{i\h}uv)}\Big)
$$
we have 
\begin{equation}
  \label{eq:lacvac}
  \begin{aligned}
(1{+}\frac{1}{m}&(z{+}\frac{uv}{i\h})){*}
e_*^{-{\frac{1}{m}(z{+}\frac{uv}{\h i})}}{*}
{\varGamma}_*(z{+}\frac{uv}{\h i})\\    
&{=}\left\{
\begin{matrix}
\prod_{k{\not=}m}(z{=}z) &
z{\not\in}-({\mathbb N}{+}m{+}\frac{1}{2})\\
\prod_{k{\not=}m}(z{=}{-}n{-}m{-}\frac{1}{2}){*}
\big(1{-} \frac{1}{n!}(\frac{1}{i\h}u)^n
{*}\varpi_{00}{*}v^{n}\big)&z{=}-(n{+}m{+}\frac{1}{2})  
\end{matrix}
\right.
  \end{aligned}
\end{equation}
where $n{\in}\mathbb N$. As opposited to the case that  
$(1{+}\frac{1}{m}(z{+}\frac{uv}{i\h}))_{+*}^{-1}{*}
\sin_*\pi(z{+}\frac{1}{i\h}uv)$ is entire function 
(cf. Proposition\,\ref{prop5.5}), there are removable 
singularities with respect to $z$.

\medskip
Multiplying 
$\prod_{k{=}1}^{\infty}
\big(1{+}\frac{1}{k}(z{+}\frac{uv}{\h i})\big)
e_*^{-{\frac{1}{k}(z{+}\frac{uv}{\h i})}}$ to both sides of 
\eqref{eq:prodinf} and using \eqref{eq:lacvac}, we have 
$$
\begin{aligned}
\lim_{N\to\infty}\prod_{k=1}^{N}&
{*}\Big(\big(1{+}\frac{1}{k}(z{+}\frac{1}{i\h}uv)\big)
{*}e_*^{{-}\frac{1}{k}(z{+}\frac{1}{i\h}uv)}\Big){*}
\varGamma_*(z{+}\frac{1}{i\h}uv)\\
&{=}
\left\{
\begin{matrix}
1 &z{\not\in}-({\mathbb N}{+}\frac{1}{2})\\
1{-}\sum_{k{=}0}^n\frac{1}{k!}
(\frac{1}{i\h}u)^k{*}\varpi_{00}{*}v^{k}, 
&z{=}-(n{+}\frac{1}{2}),   
\end{matrix}, 
\right. 
\end{aligned},
$$ 
in $H{\!o}l({\mathbb C}^2)$, where $n{\in}{\mathbb N}$.

\bigskip

\section{Products with $\sin_*\pi(z{+}\frac{1}{i\h}uv)$}

In this section we show that 
$\sin_*\pi(z{+}\frac{1}{i\h}uv){*}\varGamma_*(z{+}\frac{1}{i\h}uv)$ is
well defined as an entire function of $z$.
By recalling Euler's reflection formula, this product may be understood as 
$\frac{1}{\varGamma_*(1{-}(z{+}\frac{1}{i\h}uv))}$.
First, for ${\rm{Re}}\,z{>}-\frac{1}{2}$, we define the product 
by the integral   
\begin{equation}\label{sinsin}
\begin{aligned}
2i\sin_*\pi(z{+}\frac{1}{i\h}uv)&{*}\varGamma_*(z{+}\frac{1}{i\h}uv)\\
&{=}\lim_{T,T'\to\infty}\int_{-T}^{T'}
(e_*^{\pi i(z{+}\frac{1}{i\h}uv)}{-}e_*^{-\pi i(z{+}\frac{1}{i\h}uv)})
{*}
e^{-e^{\tau}}e_*^{\tau(z{+}\frac{uv}{i\h})}d\tau\\
&{=}
\int_{-\infty}^{\infty}
e^{-e^{\tau}}
(e_*^{(\tau{+}\pi i)(z{+}\frac{uv}{i\h})}{-}
e_*^{(\tau{-}\pi i) (z{+}\frac{uv}{i\h})})d\tau.
\end{aligned}
\end{equation}
The $\kappa$-ordering expression of \eqref{sinsin} is given as 
follows:
$$
{:}\text{\eqref{sinsin}}{:}_{\kappa}\,\,\,
{=}\,\,
\int_{{-}\infty{+}\pi i}^{\infty{+}\pi i}
e^{-e^{\tau{-}\pi i}}
e_*^{\tau(z{+}\frac{uv}{i\h})}d\tau{-}
\int_{{-}\infty{-}\pi i}^{\infty{-}\pi i}
e^{-e^{\tau{+}\pi i}}
e_*^{\tau(z{+}\frac{uv}{i\h})}d\tau.
$$
By using $e^{-e^{\tau{-}\pi i}}{=}e^{-e^{\tau{+}\pi i}}$, this 
is given by the integral 
$$
(\int_{{-}\infty{+}\pi i}^{\infty{+}\pi i}
{-}\int_{{-}\infty{-}\pi i}^{\infty{-}\pi i})
e^{e^{\tau}}
e_*^{\tau(z{+}\frac{uv}{i\h})}d\tau.
$$
Note this is not a contour integral, but it is defined for 
$\kappa{\in}{\mathbb C}{-}\{\kappa{\geq}1\}{\cup}\{\kappa{\leq}{-1}\}$. 

After this procedure, we use the analytic continuation 
via \eqref{eq:sin*}, \eqref{eq:func} to obtain the following, which 
is our main result:
\begin{thm}
  \label{main}
$\sin_*\pi(z{+}\frac{1}{i\h}uv){*}\varGamma_*(z{+}\frac{1}{i\h}uv)$
is defined as an entire function of $z$, vanishing at 
$z{\in}{\mathbb N}{+}\frac{1}{2}$ in any $\kappa$-ordering 
expression such that ${\rm{Re}}\,\kappa{<}0$, and 
$\kappa{\in}{\mathbb C}{-}\{\kappa{\geq}1\}{\cup}\{\kappa{\leq}{-1}\}$. 
\end{thm}

After careful argument about associativity,
\eqref{sinsin} can be expressed as an infinite product 
\begin{equation}
\label{prodone}
\sin_*\pi(z{+}\frac{1}{i\h}uv){*}\varGamma_*(z{+}\frac{1}{i\h}uv)
{=}
\prod_{k=1}^{\infty}{*}
\Big(\big(1{-}\frac{1}{k}(z{+}\frac{uv}{\h i})\big)_{*}
{*}e_*^{\frac{1}{k}(z{+}\frac{uv}{\h i})}\Big).
\end{equation}

Recalling the reflection formula, we may define  
$$
\frac{1}{\varGamma_*}(1{-}(z{+}\frac{1}{i\h}uv)){=} 
\sin_*\pi(z{+}\frac{1}{i\h}uv){*}\varGamma_*(z{+}\frac{1}{i\h}uv).
$$
By this we see that 
$$
\frac{1}{\varGamma}_*(1{-}(z{+}\frac{1}{i\h}uv))
\Big|_{z{=}\frac{1}{2}}=0.
$$
This supports the interpretation that 
$(\frac{1}{2}{+}\frac{1}{i\h}uv))$ is an indeterminate 
living in the set of positive integers 
${\mathbb N}{=}\{1,2,3,\cdots\}$.

\bigskip
We can form the product 
$\frac{1}{\varGamma}_*(1{-}(z{+}\frac{1}{i\h}uv)){*}
(1{-}\frac{1}{n}(z{+}\frac{1}{i\h}uv))_{-*}^{-1}
{*}e_{*}^{-\frac{1}{n}(z{+}\frac{1}{i\h}uv)}$. 
At first glance, this looks like 
$$
\prod_{k{\not=}n}{*}(1{-}\frac{1}{k}(z{+}\frac{1}{i\h}uv))
{*}e_{*}^{\frac{1}{n}(z{+}\frac{1}{i\h}uv)} 
$$
and hence as an entire function with respect to $z$. 

However, note that $(1{-}
\frac{1}{n}(z{+}\frac{1}{i\h}uv))_{-*}^{-1}$ is  
singular at $n{-}z\in {-}{\mathbb N}{-}\frac{1}{2}$, i.e. 
$z{\in}k{+}\frac{1}{2}$ for $k\geq {-}n$, and 
the same calculation as in \eqref{eq:lacvac} shows that 
$$
\sin_*\pi(z{+}\frac{1}{i\h}uv){*}
\big({\varGamma}_*(z{+}\frac{1}{i\h}uv)
{*}(1{-}\frac{1}{n}(z{+}\frac{1}{i\h}uv))_{-*}^{-1}\big)
$$ 
is not defined as an entire function, 
since some matrix elements appear in the formula as 
removable singularities. Some additional arguments is 
may be requested, since these are all removable singularities 
in the usual calculation.

\subsection{Additional support for the discrete interpretation}

We give another formula to support the discrete interpretation for 
$\frac{1}{\h i}uv$. Recall Hankel's 
\begin{picture}(130,40)(0,-10)
\thinlines
\put(0,20){\line(1,0){120}}
\put(100,0){\line(0,1){40}}
\thicklines
\put(0,22){\line(1,0){95}}
\put(0,18){\line(1,0){95}}
\put(100,20){\circle{8}}
\end{picture}
\hfill{\parbox[b]{.6\linewidth}{
formula 
$$ 
\frac{1}{\varGamma(s)}=\frac{1}{2\pi i}\int_C e^{t}t^{-s}dt, \quad
\text{(cf. \cite{ww} \,p. 244)}
$$
where $C$ is taken to be a line from 
$-\infty$ to $-\delta$, then
a circle of radius $\delta$ in the positive direction, and finally a
line from $-\delta$ to $-\infty$.}}  

Setting $s{=}\frac{1}{2}{-}\frac{1}{i\h}uv{=}
{-}\frac{1}{i\h}u{*}v$, 
we want to prove 
$\int_C e^{t}t_*^{\frac{1}{\h i}u*v}dt=0$ 
as additional support for the discrete interpretation.

\bigskip 
 By setting $t=e^{\tau{+}\pi i}$, it is easy to see that 
the Weyl ordering expression of this integral is equal to   
$$
\begin{aligned}
{:}\int_{-\infty}^0e^{t}t_*^{(\frac{1}{\h i}uv)}dt{:}_0=&
{:}\int_{-\infty}^{\infty}e^{e^{\tau{+}\pi i}}
e_*^{(\tau{+}\pi i)(1{+}\frac{1}{\h i}uv)}d\tau{:}_0\\
=&\int_{-\infty}^{\infty}e^{e^{\tau{+}\pi i}}
\frac{e^{\tau{+}\pi i}}{\cosh(\tau{+}\pi i)}
e^{\frac{1}{\h i}uv\tanh(\tau{+}\pi i)}d\tau.
\end{aligned}
$$
Hence the integral  
$$
{:}\int_{-\infty}^0e^{t}t_*^{\frac{1}{\h i}uv}dt{:}_0=
\int_{-\infty}^{\infty}e^{-e^{\tau}}
\frac{e^{\tau}}{\cosh(\tau)}
e^{\frac{1}{\h i}uv\tanh(\tau)}d\tau 
$$
exists and our integral vanishes on the axis part of $C$. 
Thus, setting $t=e^{\tau}e^{i\theta}$, we  consider 
for a fixed real $\tau <\!< 0$
$$
:A_*(\frac{uv}{\h i}):_0=
\frac{1}{2\pi}\int_0^{2\pi}
{:}e^{e^{\tau}e^{i\theta}{+}(\tau{+}i\theta)}
e_*^{(\tau{+}i\theta)(\frac{1}{\h i}uv)}{:}_0d\theta.
$$
This can be written as 
$$
\frac{1}{2\pi}\int_0^{2\pi}
e^{e^{\tau}e^{i\theta}}
\frac{e^{\tau{+}i\theta}}{\cosh(\tau{+}i\theta)}
e^{\frac{1}{\h i}uv\tanh(\tau{+}i\theta)}d\theta.
$$
We easily see the following 
\begin{lem}
\label{limit}
$\lim_{\tau\to -\infty}\frac{1}{2\pi}\int_0^{2\pi}
e^{e^{\tau}e^{i\theta}{+}(\tau{+}i\theta)}
e_*^{(\tau{+}i\theta)(\frac{2}{\h i}uv)}d\theta=0$.  
\end{lem}
Lemma\,\ref{limit} suggests that we write   
$\frac{1}{\varGamma_*}(z{+}\frac{u*v}{\h i})\Big|_{z=0}=0$, 
although this is not rigorous.

\subsection{The residue of $e_*^{t(z{+}\frac{1}{i\h}uv)}$}

We first use the Weyl ordering expression. The ${\kappa}$-ordering
expression is obtained via the intertwiner. 

\begin{lem}
\label{Resiexpct}
Let $C_k$ be a small circle of radius $\frac{\pi}{4}$ with 
the center at $\zeta=i\pi(k{+}\frac{1}{2})$. 
Then the contour integral 
$\frac{1}{2\pi i}\int_{C_k}
{:}e_*^{\zeta(z{+}\frac{1}{i\h}2uv)}{:}_0d\zeta$ gives the residue of
$e^{t}t_*^{\frac{1}{\h i}uv}$ 
and this is an entire function of 
$X=(z,\frac{1}{i\h}2uv)$. 
\end{lem}

\medskip
The continuity of the multiplication 
$(z{+}\frac{1}{i\h}2uv){*}$ requires that this function 
must satisfy the equation 
\begin{equation}
  \label{eq:muststfy}
(z{+}\frac{1}{i\h}2uv){*_0}
\int_{C_k}{:}e_*^{\zeta(z{+}\frac{1}{i\h}2uv)}{:}_0d\zeta=0,  
\end{equation}
since \eqref{eq:muststfy} equals  
$\int_{C_k}\frac{d}{d\zeta}
e_*^{\zeta(z{+}\frac{1}{i\h}2uv)}d\zeta$.
For simplicity, we set 
$$
w=\frac{1}{\h}2uv, \quad 
\Phi_k(z,w)=\frac{1}{2\pi i}\int_{C_k}
{:}e_*^{\zeta(z{+}\frac{1}{i\h}2uv)}{:}_0d\zeta.
$$
Equation \eqref{eq:muststfy} is $(iz{+}w){*_0}\Phi_k(z,w)=0$. 
Hence by the Moyal product formula, $\Phi_k(z,w)$ must satisfy 
the equation 
\begin{equation}
  \label{eq:solveq}
(iz{+}w)f(x){+}f(w)'{+}w f(w)''=0,
\end{equation}
independent of $k$. It is not difficult to see that 
equation \eqref{eq:solveq} has the unique holomorphic solution 
$f$ with initial condition $f(0)=1$.

For $f(w)=e^{aw}g(bw)$, \eqref{eq:solveq} can be rewritten as 
$$
b^2wg''(bw){+}(2abw{+}b)g'(bw){+}((a^2{+}1)w{+}a{+}iz)g(bw)=0.
$$
Thus $g(w)$ must satisfy the equation  
\begin{equation}
 \label{eq:Laguerre}
wg''(w){+}(1{+}\frac{2a}{b}w)g'(w){+}
(\frac{a^2{+}1}{b^2}w{+}\frac{a{+}iz}{b})g(w)=0.
\end{equation}

Setting $a=-\frac{1}{2}b=\pm i$, we have a Laguerre equation  
\begin{equation}
  \label{eq:LagLag}
wg''(w){+}(1{-}w)g'(w){+}\frac{1}{2}({\mp}{z}{-}1)g(w)=0,  
\end{equation}
where solution is known to be an entire function of exponential 
growth with respect to $w$.

Equation \eqref{eq:LagLag} gives two expressions for the solutions 
of \eqref{eq:solveq} using the Laguerre functions 
$L_{\nu}^{(0)}(2iw)$: 
$$
\Psi_z(w)=e^{-iw}L_{\frac{1}{2}(z{-}1)}^{(0)}(2iw), \quad
\Psi_z(w)=e^{iw}L_{-\frac{1}{2}(z{+}1)}^{(0)}(-2iw), 
$$
where 
$$
L^{(0)}_{\nu}(w)=\sum_{n=0}^{\infty}\frac{(-\nu)_n}{(n!)^2}w^n,  
\quad \nu=\frac{1}{2}({\mp}{z}{-}1).
$$
Here we use the notation
$$
(a)_n= a(a{+}1)\cdots(a{+}n{-}1), \quad (a)_0=1. 
$$

By this observation, we see that 
$\Phi_k(z,x)=c_k\Psi_z(x)$, but the
constant $c_k$ is not fixed by this method. 
To fix the constant we remark that $\Psi_z(x)$ is also
analytic in the variable $z$. 
The constant $c_k$ is fixed by investigating the case $z=0$. 

The residue of $e_*^{t\frac{1}{i\h}uv}$ is obtained in the
Weyl ordering by the contour integral 
$$
\int_{-\infty}^{\infty}(e_*^{(t{-}\pi i)\frac{1}{i\h}uv}
{-}e_*^{(t{+}\pi i)\frac{1}{i\h}uv})dt.
$$
Since $e_*^{(t{-}\pi i)\frac{1}{i\h}uv}{=}
{-}e_*^{(t{+}\pi i)\frac{1}{i\h}uv}$, 
this is given by \eqref{eq:HansenBessel}. 
\begin{lem}
\label{Residue}
The residue of 
$\frac{1}{\cosh{\zeta}}e^{(\frac{1}{i\h}\tanh{\zeta})
2uv}$  
at $\zeta=i\pi(k{+}\frac{1}{2})$ is 
$$
(-1)^k(-i)\sqrt{2\pi}J_0(\frac{2}{\h}uv),
$$
where $J_0$ is the Bessel function with the eigenvalue $0$.
\end{lem}
Comparing these we know the residue of 
$e_*^{t(z{+}\frac{1}{i\h}uv)}$.

\end{document}